\documentclass[english]{article}
\usepackage[T1]{fontenc}
\usepackage[latin9]{inputenc}
\usepackage{color}
\usepackage{amsmath}
\usepackage{amssymb}

\makeatletter
\usepackage{amsthm}\usepackage{fullpage}\usepackage{makeidx}\usepackage{amsfonts}\usepackage{latexsym}
\numberwithin{equation}{section}

\numberwithin{figure}{section}

\usepackage{fullpage}\usepackage{makeidx}\usepackage{amsfonts}\usepackage{latexsym}

\def\Z{{\mathbb{Z}}}

\theoremstyle{definition}
\newtheorem{lemma}{Lemma}[section]
\newtheorem{theorem}[lemma]{Theorem}\newtheorem{proposition}[lemma]{Proposition}\newtheorem{definition}[lemma]{Definition}\newtheorem{remark}[lemma]{Remark}\usepackage{times}

\usepackage{babel}

\usepackage[blocks]{authblk}
\title{$S$-matrix in permutation orbifolds}

\author{Chongying Dong\footnote{Supported by the Simons Foundation 634104}}
\affil{Department of Mathematics, University of
California, Santa Cruz, CA 95064 USA}

\author{Feng Xu }

\affil{University of California at Riverside, Riverside, CA {\rm 92521} USA }

\author{Nina Yu\footnote{Supported by National Natural Science Foundation of China  11971396 and  12131018}}
\affil{School of Mathematical Sciences, Xiamen University, Xiamen, Fujian 361005, CHINA}

\numberwithin{equation}{section}

\usepackage{babel}

\usepackage{babel}

\makeatother

\usepackage{babel}
\begin{document}
\maketitle
\begin{abstract}
For a fixed positive integer $k$, any element $g$ of the permutation
group $S_{k}$ acts on the tensor product vertex operator algebra
$V^{\otimes k}$ in the obvious way. In this paper, we determine the
$S$-matrix of $\left(V^{\otimes k}\right)^{G}$ if $G=\left\langle g\right\rangle $
is the cyclic group generated by $g=\left(1,\ 2,\cdots,k\right).$
\end{abstract}

\section{Introduction}

Let $V$ be a vertex operator algebra. Permutation orbifold theory studies the representations
of the tensor product vertex operator algebra $V^{\otimes k}$ with
the natural action of the symmetric group $S_{k}$ as an automorphism
group, where $k$ is a positive integer. In this paper we determine the
$S$-matrix of $\left(V^{\otimes k}\right)^{G}$ where $G=\left\langle g\right\rangle $
is the cyclic group generated by $g=\left(1,\ 2,\cdots,k\right).$

The study of permutation orbifolds was initiated in \cite{BHS}, where the twisted modules,
genus one characters and the fusion rules for cyclic permutations
for affine vertex operator algebras and the Virasoro vertex operator
algebras were studied. The genus one characters and modular transformation
properties of permutation orbifolds for a general rational conformal
field theory were given in \cite{Ba}. The twisted
modules for $V^{\otimes k}$
were constructed for any  permutation automorphism of  $V^{\otimes k}$
in \cite{BDM}. Specifically, let $g$ be a $k$-cycle, which is naturally
an automorphism of $V^{\otimes k},$ then for any $V$-module $\left(W,Y_{W}\left(\cdot,z\right)\right),$
a canonical $g$-twisted $V^{\otimes k}$-module structure on $W$
was obtained. Furthermore, it was proved that there is an isomorphism
of the categories of weak, admissible and ordinary $V$-modules and
the categories of weak, admissible and ordinary $g$-twisted $V^{\otimes k}$-modules, respectively. The $C_{2}$-cofiniteness of permutation orbifolds and
general cyclic orbifolds was established later in \cite{A1,A2,M1,M2}.
An equivalence of two constructions \cite{FLM,Le,BDM} of twisted
modules for permutation orbifolds of lattice vertex operator algebras
was given in \cite{BHL}. The permutation orbifolds of the lattice
vertex operator algebras with $k=2$ and $k=3$ were extensively studied
in \cite{DXY1,DXY2,DXY3}. Fusion products of $V^{\otimes k}$-modules
with $\sigma$-twisted $V^{\otimes k}$-module for any $\sigma\in S_{k}$
were studied in \cite{DLXY}.

It is well known that the modular group $\Gamma=SL_{2}\left(\mathbb{Z}\right)$ acts on the conformal block
of a rational, $C_2$-cofinite vertex operator algebra \cite{Z, DLM00, DLN}. The action of $S=\left(\begin{array}{cc}
0 & -1\\
1 & 0
\end{array}\right)$ is called the {\em $S$-matrix}, which is key
to understand the action of $\Gamma$.  Since the conformal block of $V^{G}$ spanned by the
trace functions on the irreducible $V^{G}$-modules appearing in twisted
$V$-modules is equal to the twisted conformal block of $V$ spanned
by the trace functions on the irreducible twisted modules \cite{DRX17}, one can give a precise formula for the restricted $S$-matrix in terms of the $S$-matrix of the twisted $V$-modules, where the \emph{restricted $S$-matrix} of
$V^{G}$ is the  restriction of
the $S$-matrix of $V^{G}$ to the irreducible $V^{G}$-modules appearing
in the twisted modules \cite{DRX21}. Some entries
of the restricted $S$-matrix have been computed in \cite{DRX17}
for studying quantum dimensions and global dimensions for vertex operator
algebras $V^{G}$.
As pointed out in \cite{DRX21}, the restricted $S$-matrix is equal to  the $S$-matrix of $V^G$ if $V^G$ is rational and $C_2$-cofinite. The $S$-matrix for cyclic group $G$ and holomorphic vertex operator algebra has been studied in \cite{EMS} for constructing holomorphic vertex operator algebra with central charge 24.

From \cite{CM}, if $V$ is rational and $C_2$-cofinite, $G$ is an abelian automorphism group of $V$ then $V^G$ is also rational and $C_2$-cofinite. Note that $V^{\otimes k}$ is also rational and $C_2$-cofinite. This implies that  $\left(V^{\otimes k}\right)^{G}$  is rational and $C_2$-cofinite if $G$ is cyclic. Our main result in this paper is an explicit formula of $S$-matrix of  $\left(V^{\otimes k}\right)^{G}$ in terms of the action of $\Gamma$ on the conformal block of $V$ if $G$ is generated by $g=(1,2,...,k).$ If $k=2$ or $k$ is a general prime, this result has been obtained previously in \cite{BHS} and \cite{DRX21}, respectively.  The main idea is that the action of $\Gamma$ on the twisted conformal block for $V^{\otimes k}$ is determined explicitly in terms of the action of $\Gamma$ on the conformal block of $V$
as the twisted modules for $V^{\otimes k}$ are known by using the $V$-modules \cite{BDM}.
If $k$ is not a prime, some powers of $g$ may not be $k$-cycles. This makes the computation of $S$-matrix of $\left(V^{\otimes k}\right)^{G}$ much more complicated.

This paper is organized as follows. We present basic notions
and results on vertex operator algebras in Section 2. Twisted modules  in
permutation orbifolds are discussed in Section 3. The trace functions of twisted modules
for $V^{\otimes k}$-modules are computed  and the $S$-matrix on twisted conformal block
for $V^{\otimes k}$ is given in Section 4. In Section 5, we give an explicit
expression for  $S$-matrix of $\left(V^{\otimes k}\right)^{\left\langle g\right\rangle }$
in terms of the $S$-matrix of $V^{\otimes k}$ obtained in Section
4.

\section{Preliminary}

In this section, we review the basics on vertex operators algebras.

\subsection{Basics}

Let $V=(V,\ Y,\ \mathbf{1},\ \omega)$ be a vertex operator algebra.
Let $Y(v,\ z)=\sum_{n\in\mathbb{Z}}v_{n}z^{-n-1}$ denote the vertex
operator of $v\in V$, where $v_{n}\in\mbox{End}(V)$. We first recall
some basic notions from \cite{FLM,Z,DLM96,DLM97}.

\begin{definition} An \emph{automorphism} $g$ of a vertex operator
algebra $V$ is a linear isomorphism of $V$ satisfying $g\left(\omega\right)=\omega$
and $gY\left(v,z\right)g^{-1}=Y\left(gv,z\right)$ for any $v\in V$.
We denote by $\mbox{Aut}\left(V\right)$ the group of all automorphisms
of $V$. \end{definition}

For a subgroup $G\le\mbox{Aut}\left(V\right)$ the fixed point set
$V^{G}=\left\{ v\in V\mid g\left(v\right)=v,\forall g\in G\right\} $
has a vertex operator algebra structure. Let $g$ be an automorphism
of a vertex operator algebra $V$ of order $T$. Denote the decomposition
of $V$ into eigenspaces of $g$ as
\[
V=\oplus_{r\in\mathbb{Z}/T\text{\ensuremath{\mathbb{Z}}}}V^{r}
\]
where $V^{r}=\left\{ v\in V\mid gv=e^{-2\pi ir/T}v\right\} $.

\begin{definition} A \emph{weak $g$-twisted $V$-module} $M$ is
a vector space with a linear map
\begin{align*}
Y_{M}:  & V\to\left(\text{End}M\right)\{z\}\\
 & v\mapsto Y_{M}\left(v,z\right)=\sum_{n\in\mathbb{Q}}v_{n}z^{-n-1}\ \left(\text{where}\  v_{n}\in\mbox{End}M\right)
\end{align*}
which satisfies the following: for all $0\le r\le T-1$, $u\in V^{r}$,
$v\in V$, $w\in M$,
\[
Y_{M}\left(u,z\right)=\sum_{n\in-\frac{r}{T}+\mathbb{Z}}u_{n}z^{-n-1},
\]
\[
u_{l}w=0\ \text{for}\ l\gg0,
\]
\[
Y_{M}\left(\mathbf{1},z\right)=\text{Id}_{M},
\]
\[
z_{0}^{-1}\text{\ensuremath{\delta}}\left(\frac{z_{1}-z_{2}}{z_{0}}\right)Y_{M}\left(u,z_{1}\right)Y_{M}\left(v,z_{2}\right)-z_{0}^{-1}\delta\left(\frac{z_{2}-z_{1}}{-z_{0}}\right)Y_{M}\left(v,z_{2}\right)Y_{M}\left(u,z_{1}\right)
\]

\begin{equation}
=z_{2}^{-1}\left(\frac{z_{1}-z_{0}}{z_{2}}\right)^{-r/T}\delta\left(\frac{z_{1}-z_{0}}{z_{2}}\right)Y_{M}\left(Y\left(u,z_{0}\right)v,z_{2}\right),\label{Jacobi for twisted V-module}
\end{equation}
where $\delta\left(z\right)=\sum_{n\in\mathbb{Z}}z^{n}$.

\end{definition}

\begin{definition}

An \emph{admissible $g$-twisted $V$-module} is a $\frac{1}{T}\mathbb{Z}_{+}$-graded
weak $g$-twisted $V$-module $M$: $M=\oplus_{n\in\frac{1}{T}\mathbb{Z}_{+}}M\left(n\right)$
such that $v_{m}M\left(n\right)\subseteq M\left(n+\text{wt}v-m-1\right)$
for homogeneous $v\in V$ and $m,n\in\frac{1}{T}\mathbb{Z}$.

\end{definition}

If $M=\oplus_{n\in\frac{1}{T}\mathbb{Z}_{+}}M\left(n\right)$ is an
admissible $g$-twisted $V$-module, the contragredient module $M'$
is defined as follows:
\[
M'=\oplus_{n\in\frac{1}{T}\mathbb{Z}_{+}}M\left(n\right)^{\ast},
\]
where $M\left(n\right)^{\ast}=\text{Hom}_{\mathbb{C}}\left(M\left(n\right),\mathbb{C}\right).$
The vertex operator $Y_{M'}\left(a,z\right)$ is defined for $a\in V$
via

\[
\left\langle Y_{M'}\left(a,z\right)f,u\right\rangle =\left\langle f,Y_{M}\left(e^{zL\left(1\right)}\left(-z^{-2}\right)^{L\left(0\right)}a,z^{-1}\right)u\right\rangle ,
\]
where $\left\langle f,w\right\rangle =f\left(w\right)$ is the natural
pairing $M'\times M\to\mathbb{C}$. One can prove the following \cite{FHL, X}:

\begin{lemma} $\left(M',Y_{M'}\right)$ is an admissible $g^{-1}$-twisted
$V$-module.

\end{lemma}

\begin{definition}

A $g$-\emph{twisted $V$-module} is a weak $g$-twisted $V$-module
$M$ which carries a $\mathbb{C}$-grading induced by the spectrum
of $L(0)$ where $L(0)$ is the component operator of $Y(\omega,z)=\sum_{n\in\mathbb{Z}}L(n)z^{-n-2}.$
That is, we have $M=\bigoplus_{\lambda\in\mathbb{C}}M_{\lambda},$
where $M_{\lambda}=\left\{ w\in M\mid L(0)w=\lambda w\right\} $. Moreover,
$\dim M_{\lambda}$ is finite and for fixed $\lambda,$ $M_{\frac{n}{T}+\lambda}=0$
for all small enough integers $n.$

\end{definition}

If $g=\text{Id}_{V}$ we have the notions of weak, admissible and ordinary
$V$-modules \cite{DLM98}.

\begin{definition} A vertex operator algebra $V$ is said to be \emph{regular}
if the weak $V$-module category is semisimple.

\end{definition}

\begin{definition} A vertex operator algebra $V$ is said to be \emph{$g$-rational}
if the admissible $g$-twisted module category is semisimple. We say
$V$ is \emph{ rational} if $V$ is $1$-rational.

\end{definition}

\begin{definition} A vertex operator algebra $V$ is said to be \emph{$C_{2}$-cofinite}
if $V/C_{2}\left(V\right)$ is finite dimensional, where $C_{2}\left(V\right)=\left\langle u_{-2}v\mid u,v\in V\right\rangle $.

\end{definition}

\begin{definition} A vertex operator algebra $V=\oplus_{n\in\mathbb{Z}}V_{n}$
is said to be of \emph{CFT type} if $V_{n}=0$ for negative $n$ and $V_{0}=\mathbb{C}\boldsymbol{1}.$

\end{definition}

If $M=\oplus_{n\in\frac{1}{T}\mathbb{Z}_{+}}M(n)$ is an irreducible
admissible $g$-twisted $V$-module, then there is a complex number
$\lambda_{M}$ such that $L(0)|_{M(n)}=\lambda_{M}+n$ for all $n.$
As a convention, we assume $M(0)\ne0$, and $\lambda_{M}$ is called
the {\em weight} or {\em conformal weight} of $M.$

\begin{remark} (1) If $V$ is rational then there are only finitely
irreducible admissible $V$-modules up to isomorphism and each irreducible
admissible $V$-module is ordinary \cite{DLM98}.

(2) If $V$ is of CFT type, then regularity is equivalent to rationality
and $C_{2}$-cofiniteness \cite{KL,ABD}.

(3) Assume that $V$ is rational and $C_{2}$-cofinite. Then $V$
is $g$-rational for any finite automorphism $g$ \cite{ADJR}, and
it was proved in \cite{DLM00} that $\lambda_{M}$ is a rational number
for every irreducible $g$-twisted $V$-module $M$.

\end{remark}

In the rest of this paper, we assume that $V=\oplus_{n\ge0}V_{n}$ is a simple, rational, $C_{2}$-cofinite
vertex operator algebra of CFT type and  $G$ is  a finite automorphism
group of  $V$ such that the conformal weight of any irreducible $g$-twisted $V$-module
$M$ is nonnegative and is zero if and only if $M=V$. Under the above assumptions, $V^{G}$ is rational and $C_{2}$-cofinite if $G$ is solvable \cite{CM,M2}.

\subsection{Modular invariance \label{modular invariance}}

Now we review some results on modular invariance in orbifold theory from \cite{Z,DLM00}. These results play important roles in this paper.

We need  the action of Aut$\left(V\right)$ on the set of twisted modules.
Let $g,h\in\text{Aut}\left(V\right)$ with $g$ finite order. If $\left(M,Y_{M}\right)$
is a weak $g$-twisted $V$-module, there is a weak $h^{-1}gh$-twisted
$V$-module $\left(M\circ h,Y_{M\circ h}\right)$ where $M\circ h\cong M$
as vector spaces and $Y_{M\circ h}\left(v,z\right)=Y_{M}\left(hv,z\right)$
for $v\in V$. This defines a right action of $\text{Aut}\left(V\right)$
on the set of weak twisted $V$-modules and on isomorphism classes of weak twisted
$V$-modules.
$M$ is called \emph{$h$-stable} if $M$ and $M\circ h$ are isomorphic.

Assume that $g,h$ commute. Then $h$ acts on the $g$-twisted modules.
Denote by $\mathcal{\text{\ensuremath{\mathcal{\text{\ensuremath{\mathfrak{U}\left(g\right)}}}}}}$
the equivalence classes of irreducible $g$-twisted $V$-modules and
\[\mathfrak{U}\left(g,h\right)=\left\{ M\in\mathfrak{U}\left(g\right)\mid M\circ h\cong M\right\} .\]
Both $\mathfrak{U}\left(g\right)$ and $\mathfrak{U}\left(g,h\right)$
are finite sets since $V$ is $g$-rational for all $g$.

Let $M$ be an irreducible $g$-twisted $V$-module and $G_{M}$ be
a subgroup of $G$ consisting of $h\in G$ such that $M\circ h$ and
$M$ are isomorphic. By Schur's Lemma there is a projective representation
$\phi$ of $G_{M}$ on $M$ such that
\[
\phi\left(h\right)Y\left(u,z\right)\phi\left(h\right)^{-1}=Y\left(hu,z\right)
\]
for $h\in G_{M}$. If $h=1$ we take $\phi\left(1\right)=1$. Note
that $g$ lies in $G_{M}$ as $g$ acts naturally on any admissible $g$-twisted-module
$M$ such that $g|_{M(n)}=e^{2\pi in}$ for $n\in\frac{1}{T}\Z.$
We will use this action of $g$ throughout this paper.

Set $o\left(v\right)=v_{\text{wt}v-1}$ for homogeneous $v\in V$.
Then $o\left(v\right)$ is a degree zero operator of $v$. Let $\mathbb{H}$
be the complex upper half-plane. Here and below we set $q=e^{2\pi i\tau}$
where $\tau\in\mathbb{H}$. For $v\in V$, set

\begin{equation}
Z_{M}\left(v,\left(g,h\right),\tau\right)=\text{tr}_{M}o\left(v\right)\phi\left(h\right)q^{L\left(0\right)-c/24}=q^{\lambda-c/24}\sum_{n\in\frac{1}{T}\mathbb{Z}_{+}}\text{tr}_{M_{\lambda+n}}o\left(v\right)\phi\left(h\right)q^{n}.\label{def of trace function}
\end{equation}
Then $Z_{M}\left(v,\left(g,h\right),\tau\right)$ is a holomorphic
function on $\mathbb{H}$ \cite{Z, DLM00}. We write $Z_{M}\left(v,\tau\right)=Z_{M}\left(v,\left(g,1\right),\tau\right)$
for short. Then $\chi_{M}\left(\tau\right)=Z_{M}\left(\boldsymbol{1},\tau\right)$
is called the \emph{character }of $M$.

Recall that there is another vertex operator algebra $\left(V,Y\left[\ \ \ \right],\boldsymbol{1},\tilde{\omega}\right)$
associated to $V$ (see \cite{Z}). Here $\tilde{\omega}=\omega-c/24$
and for homogeneous $v\in V$,
\[
Y\left[v,z\right]=Y\left(v,e^{z}-1\right)e^{z\cdot\text{wt}v}=\sum_{n\in\mathbb{Z}}v\left[n\right]z^{n-1}.
\]
We write
\[
Y\left[\tilde{\omega},z\right]=\sum_{n\in\mathbb{Z}}L\left(n\right)z^{-n-2}.
\]
The weight of a homogeneous $v\in V$ in the second vertex operator
algebra is denoted by $\text{wt}\left[v\right].$

The modular group $\Gamma=SL_{2}\left(\mathbb{Z}\right)$ is the group
of $2\times2$ integral matrices with determinant 1. Denote by $\Gamma\left(N\right)$
the kernel of the reduction modulo $N$ epimorphism $\pi_{N}:SL_{2}\left(\mathbb{Z}\right)\to SL_{2}\left(\mathbb{Z}_{N}\right)$.
A subgroup $G_{N}$ of $SL_{2}\left(\mathbb{Z}\right)$ is called
a \emph{congruence subgroup of level N }if $N$ is the least positive
integer such that $\Gamma\left(N\right)\le G_{N}$.

Let $P\left(G\right)$ be the set of the ordered commutating pairs
in $G$. For $\left(g,h\right)\in P\left(G\right)$ and $M\in\mathfrak{U}$$\left(g,h\right),$
$Z_{M}\left(v,\left(g,h\right),\tau\right)$ is a function on $V\times\mathbb{H}.$
Let $W$ be the vector space spanned by such functions. Then by \cite{DLM00}
the dimension of $W$ is equal to $\sum_{\left(g,h\right)\in P\left(G\right)}\left|\mathfrak{U}\left(g,h\right)\right|$.
Now we define an action of the modular group $\Gamma$ on $W$ such
that
\[
Z_{M}|_{\gamma}\left(v,\left(g,h\right),\tau\right)=\left(c\tau+d\right)^{-\text{wt}\left[v\right]}Z_{M}\left(v,\left(g,h\right),\gamma\tau\right),
\]
where $\gamma:\tau\,\mapsto\frac{a\tau+b}{c\tau+d},$ $\gamma=\left(\begin{array}{cc}
a & b\\
c & d
\end{array}\right)\in\Gamma=SL\left(2,\mathbb{Z}\right).$ Let $\gamma\in\Gamma$ act on the right of $P\left(G\right)$ via
\[
\left(g,h\right)\gamma=\left(g^{a}h^{c},g^{b}h^{d}\right).
\]

We will need to use the following results from \cite{DLM00,Z,DLN,DR}:

\begin{theorem}\label{modular invariance thm} Let $V$, $G$ and $W$
be as before. Then

(1) There is a representation $\rho:\Gamma\to GL\left(W\right)$ such
that for $\left(g,h\right)\in P\left(G\right)$, $\gamma=\left(\begin{array}{cc}
a & b\\
c & d
\end{array}\right)\in\Gamma$ and $M\in\mathfrak{U}\left(g,h\right),$
\[
Z_{M}|_{\gamma}\left(v,\left(g,h\right),\tau\right)=\sum_{N\in\mathfrak{U}\left(g^{a}h^{c},g^{b}h^{d}\right)}\gamma_{M,N}Z_{N}\left(v,\left(g,h\right),\tau\right),
\]
where $\rho\left(\gamma\right)=\left(\gamma_{M,N}\right)$. That is,

\[
Z_{M}\left(v,\left(g,h\right),\gamma\tau\right)=\left(c\tau+d\right)^{\text{wt}\left[v\right]}\sum_{N\in\mathfrak{U}\left(g^{a}h^{c},g^{b}h^{d}\right)}\gamma_{M,N}Z_{N}\left(v,\left(g^{a}h^{c},g^{b}h^{d}\right),\tau\right).
\]

(2) The cardinalities $\left|\mathfrak{U}\left(g,h\right)\right|$
and $\left|\mathfrak{U}\left(g^{a}h^{c},g^{b}h^{d}\right)\right|$
are equal for any $\left(g,h\right)\in P\left(G\right)$ and $\gamma\in\Gamma$.
In particular, the number of irreducible $g$-twisted $V$-modules
exactly equals the number of irreducible $V$-modules that are $g$-stable.

(3) Each $Z_{M}\left(v,\left(g,h\right),\tau\right)$ is a modular
form of weight $\text{wt}\left[v\right]$ on the congruence subgroup.
In particular, the character $\chi_{M}\left(\tau\right)$ is a modular
function on the same congruence subgroup.

\end{theorem}

Since the modular group $\Gamma$ is generated by $S=\left(\begin{array}{cc}
0 & -1\\
1 & 0
\end{array}\right)$ and $T=\left(\begin{array}{cc}
1 & 1\\
0 & 1
\end{array}\right)$, the representation $\rho$ is uniquely determined by $\rho\left(S\right)$
and $\rho\left(T\right).$ The matrix $\rho\left(S\right)$ is called
the\emph{ $S$-matrix } of the orbifold theory. Consider a special
case of the $S$-transformation:

\[
Z_{M}\left(v,-\frac{1}{\tau}\right)=\tau^{\text{wt}\left[v\right]}\sum_{N\in\mathfrak{U}\left(1,g^{-1}\right)}S_{M, N}Z_{N}\left(v,\left(1,g^{-1}\right),\tau\right)
\]
for $M\in\mathfrak{U}\left(g\right)$ and

\[
Z_{N}\left(v,\left(1,g\right),-\frac{1}{\tau}\right)=\tau^{\text{wt}\left[v\right]}\sum_{M\in\mathfrak{U}\left(g\right)}S_{N, M}Z_{M}\left(v,\tau\right)
\]
for $N\in\mathfrak{U}\left(1\right)$. The matrix $S=\left(S_{M,N}\right)_{M,N\in\mathfrak{U}\left(1\right)}$ is
called the \emph{$S$-matrix of  $V$. }

\section{Twisted modules in permutation orbifold}

In the rest of this paper, we fix $g=(1,\ 2,\cdots,k)$ where $k$
is a positive integer. Now we study twisted modules of the tensor
product vertex operator algebra $V^{\otimes k}$ under cyclic permutation
group $\left\langle g\right\rangle $. In this section, we first review
the structure of $g$-twisted $V^{\otimes k}$-modules from \cite{BDM}.
Then we study $g^{s}$-twisted $V^{\otimes k}$-modules, for $2\le s<k$.
Note that here $g^{s}$ could be a product of several disjoint cycles.

\subsection{\label{g-twisted module for tensor product VOA}$g$-twisted $V^{\otimes k}$-modules}

It is proved in \cite{DLM00} that the number of irreducible $V$-modules
is equal to the number of irreducible $g$-twisted $V^{\otimes k}$-modules
up to isomorphism. A functor $T_{g}$ from the category of $V$-modules
to the category of $g$-twisted $V^{\otimes k}$-modules is constructed
in \cite{BDM}. Recall that
\[
\Delta_{k}\left(z\right)=\exp\left(\sum_{n\ge1}a_{n}z^{-n/k}L\left(n\right)\right)k^{-L\left(0\right)}z^{\left(1/k-1\right)L\left(0\right)},
\]
where the coefficients $a_{n}$ for $n\ge1$ are uniquely determined
by
\[
\exp\left(\sum_{n\ge1}-a_{n}x^{n+1}\frac{d}{dx}\right)x=\frac{1}{k}\left(1+x\right)^{k}-\frac{1}{k}.
\]

For $v\in V$ we denote by $v^{j}\in V^{\otimes k}$ the vector whose
$j$-th tensor factor is $v$ and whose other tensor factors are 1:
\[
v^{j}=1^{\otimes\left(j-1\right)}\otimes v\otimes1^{\otimes\left(k-j\right)}.
\]
We have $Y\left(v^{j},z\right)=1^{\otimes\left(j-1\right)}\otimes Y\left(v,z\right)\otimes1^{\otimes\left(k-j\right)}$.
Note that $gv^{j}=v^{j+1}$ for $j=1,\cdots,k,$ where $v^{k+1}=v^{1}$
by convention.

For any $V$-module $\left(W,Y_{W}\right)$, there is a $g$-twisted
$V^{\otimes k}$-module $\left(T_{g}^{k}\left(W\right),Y_{T_{g}^{k}\left(W\right)}\right)$,
where $T_{g}^{k}\left(W\right)=W$ as a vector space and the vertex
operator map $Y_{T_{g}^{k}\left(W\right)}\left(\cdot,z\right)$ is
uniquely determined by

\[
Y_{T_{g}^{k}\left(W\right)}\left(u^{1},z\right)=Y_{W}\left(\Delta_{k}\left(z\right)u,z^{1/k}\right)\ \text{for}\ u\in V.
\]
Furthermore, every $g$-twisted $V^{\otimes k}$-module is isomorphic
to one of this form.  Assume that $W^0, W^1, \cdots,W^{p}$ are all
the irreducible $V$-modules. Now we see that $T_{g}^{k}\left(W^{0}\right),$
$T_{g}^{k}\left(W^{1}\right),\cdots,$ $T_{g}^{k}\left(W^{p}\right)$
are all the irreducible $g$-twisted $V^{\otimes k}$-modules. There
is a $\frac{1}{k}\mathbb{Z}_{+}$-gradation on $T_{g}^{k}\left(W^{i}\right)$
such that $T_{g}^{k}\left(W^{i}\right)=\oplus_{n\ge0}T_{g}^{k}\left(W^{i}\right)\left(\frac{n}{k}\right)$
with $T_{g}^{k}\left(W^{i}\right)\left(\frac{n}{k}\right)\cong W^{i}\left(n\right)$
as a vector space, and $Y_{g}\left(v,z\right)=\sum_{m\in\frac{1}{k}\mathbb{Z}}v_{m}z^{-m-1}$
for $v\in V^{\otimes k}$.

\subsection{$g^{s}$-twisted $V^{\otimes k}$-modules}

For any $1\le s<k$, we now consider $g^{s}$-twisted $V^{\otimes k}$-module.
Let $d=\text{gcd}\left(s,k\right)$ where $\text{gcd}\left(s,k\right)$
is the greatest common divisor of $s$ and $k$. Let $m=\frac{s}{d}$,
$l=\frac{k}{d}$. Then $\text{gcd}\left(m,l\right)=1$ and $o\left(g^{s}\right)=l$.

\subsubsection{Structure of $g^{s}$-twisted modules}

Note that

\[
g^{d}=h_{1}\cdots h_{d},
\]
where $h_i, 1\le i\le d$ are $l$-cycles:
\begin{gather*}
h_{1}=\left(1,d+1,2d+1,\cdots,\left(l-1\right)d+1\right),\\
h_{2}=\left(2,d+2,2d+2,\cdots,\left(l-1\right)d+1\right),\\
\vdots\\
h_{d}=\left(d,2d,\cdots,ld\right).
\end{gather*}
Therefore
\begin{equation}
g^{s}=g^{dm}=h_{1}^{m}h_{2}^{m}\cdots h_{d}^{m}\label{cycle decompositon of g^s},
\end{equation}
where each $h_{i}^{m}$ is\textcolor{magenta}{{} }an $l$-cycle as
$\text{gcd}\left(m,l\right)=1$. From now on, we let $\text{Irr}(V)=\{W^0, W^1,\cdots, W^p\}$ be the set of all irreducible $V$-modules. Then by \cite{BDM}, each irreducible $g^{s}$-twisted
module is of the form
\[
T_{h_{1}^{m}}^{l}\left(M^{1}\right)\otimes\cdots\otimes T_{h_{d}^{m}}^{l}\left(M^{d}\right),
\]
where $M^{1},\cdots,M^{d}\in\text{Irr}\left(V\right)$.
Note that each $T_{h_{i}^{m}}^{l}\left(M^{i}\right)$ is a $h_{i}^{m}$-twisted
$V^{\otimes l}$-module. For short we will write
\begin{equation}
T_{g^{s}}^{M^{1},\cdots,M^{d}}=T_{h_{1}^{m}}^{l}\left(M^{1}\right)\otimes\cdots\otimes T_{h_{d}^{m}}^{l}\left(M^{d}\right).\label{g^s-module general form}
\end{equation}
Since $gh_{i}g^{-1}=h_{i+1}$ or $g^{-1}h_{i+1}g=h_{i}$, we obtain
\[
g^{-1}g^{d}g=\left(g^{-1}h_{1}g\right)\left(g^{-1}h_{2}g\right)\cdots\left(g^{-1}h_{d}g\right)=h_{d}h_{1}\cdots h_{d-1}.
\]
Since $s=dm,$ we get
\[
g^{-1}g^{s}g=\left(g^{-1}g^{d}g\right)^{m}=h_{d}^{m}h_{1}^{m}\cdots h_{d-1}^{m}.
\]
By \cite{BDM} we see that
\[
T_{g^{s}}^{M^{1},\cdots,M^{d}}\circ g\cong T_{g^{s}}^{M^{d},M^{1},\cdots,M^{d-1}}.
\]
In particular,
\begin{equation}
T_{g^{s}}^{M^{1},\cdots, M^{d}}\circ g^{d}\cong T_{g^{s}}^{M^{1},\cdots, M^{d}}.\label{g^d stabilizer}
\end{equation}
This produces a $d$-cycle $\sigma=\left(1,d,d-1,\cdots,2\right).$
Let $1\le r<k$ and $f=\text{gcd}$$\left(d,r\right)$. Set $a=\frac{r}{f}$
and $b=\frac{d}{f}$. Then $\text{gcd}\left(a,b\right)=1$ and $o\left(\sigma^{r}\right)=b$. By similar argument as above, we have
\[
\sigma^{r}=\sigma^{af}=\sigma_{1}^{a}\cdots\sigma_{f}^{a},
\]
where $\sigma_i, 1\le i\le f$ are $b$-cycles:
\begin{gather*}
\sigma_{1}=\left(1,f+1,2f+1,\cdots,\left(b-1\right)f+1\right),\\
\sigma_{2}=\left(2,f+2,2f+2,\cdots,\left(b-1\right)f+2\right),\\
\vdots\\
\sigma_{f}=\left(f,2f,3f,\cdots,bf\right).
\end{gather*}
This implies that $T_{g^{s}}^{M^{1},\cdots,M^{d}}\circ g^{r}\cong T_{g^{s}}^{M^{1},\cdots,M^{d}}$
if and only if $M^{i}=M^{i+jf}$for $i=1,...,f$ and $j=1,...,b-1.$
Hence we have the following lemma.

\begin{lemma} \label{g^s, g^r modules} Let $1\le r,s<k$, $d=\text{gcd}\left(s,k\right)$,
$f=\text{gcd}$$\left(d,r\right)$. Set $l=\frac{k}{d},$ $m=\frac{s}{d}$
and $b=\frac{d}{f}$. Then the equivalence classes of irreducible
$g^{s}$-twisted $V^{\otimes k}$-modules that are $g^{r}$-stable
are
\begin{equation}
\mathfrak{U}\left(g^{s},g^{r}\right)=\left\{ \left(T_{h_{1}^{m}}^{l}\left(M^{1}\right)\right)^{\otimes b}\otimes\cdots\otimes\left(T_{h_{f}^{m}}^{l}\left(M^{f}\right)\right)^{\otimes b}\mid M^{1},\cdots,M^{f}\in\text{Irr}\left(V\right)\right\} ,\label{general module g^s,g^r}
\end{equation}
where each $h_{i}^{m}$ is an $l$-cycle, $1\le i\le f$.

\end{lemma}

\subsubsection{Weights of $g^{s}$-twisted $V^{\otimes k}$-modules }

Let $k,d,l$ be as before. Set $\overline{v_{i}}=\sum_{j=1}^{l}v_{i}^{j}$
where $v_{i}^{j}\in V^{\otimes l}$ denotes the vector whose $j$-th
tensor factor is $v_{i}\in V$ and whose other tensor factors are
1. (The subindex $i$ is to distinguish tensor factor in tensor product
of vectors of such form, as we will see later.) Then $\bar{v}=\overline{v_{1}}\otimes1^{\left(d-1\right)l}+1^{\otimes l}\otimes\overline{v_{2}}\otimes1^{\otimes\left(d-2\right)l}+\cdots+1^{\left(d-1\right)l}\otimes\overline{v_{d}}$.
Let $\omega$ be the Virasoro vector of $V$. We have
$\Delta_{l}\left(z\right)\omega=\frac{z^{2\left(1/l-1\right)}}{l^{2}}\left(\omega+\frac{\left(l^{2}-1\right)c}{24}z^{-2/l}\right)$ \cite{BDM}.
Write $Y_{h_{i}^{m}}\left(\text{\ensuremath{\overline{\omega_{i}}}},z\right)=\sum_{n\in\mathbb{Z}}L_{h_{i}^{m}}\left(n\right)z^{-n-2}.$
Then $L_{h_{i}^{m}}\left(0\right)=\frac{1}{l}L\left(0\right)+\frac{\left(l^{2}-1\right)c}{24l}$.

Now we compute the weight of the $g^{s}$-twisted
module $\mathcal{M}=T_{g^{s}}^{M^{1},\cdots,M^{d}}$.  Let $Y_{\mathcal{M}}\left(\bar{\omega},z\right)=\sum_{n\in\mathbb{Z}}L_{\mathcal{M}}\left(n\right)z^{-n-2},$
where $\bar{\omega}=\sum_{j=1}^{k}\omega^{j}$ is the Virasoro vector
of $V^{\otimes k}.$ It is easy to see  that
\begin{equation}
L_{\mathcal{M}}\left(0\right)=\sum_{i=1}^{d}L_{h_{i}^{m}}\left(0\right)=\sum_{i=1}^{d}\frac{1}{l}L^{i}\left(0\right)+\frac{d\left(l^{2}-1\right)c}{24l},\label{virasoro of M}
\end{equation}
where $L^{i}(0)$ is the $L(0)$ on the $i$-tensor factor of ${\cal M}.$

Denote the weight of $M^{i}$ by $\lambda_{i}$ and the weight of
$\mathcal{M}$ by $\lambda_{\mathcal{M}}$. Then
\begin{alignat}{1}
\lambda_{\mathcal{M}}=\frac{\lambda_{1}}{l}+\frac{\left(l^{2}-1\right)c}{24l}+\cdots+\frac{\lambda_{d}}{l}+\frac{\left(l^{2}-1\right)c}{24l} & =\frac{\lambda_{1}+\cdots+\lambda_{d}}{l}+\frac{d\left(l^{2}-1\right)c}{24l}\label{weight of M}
\end{alignat}
and
\[
\mathcal{M}=\bigoplus_{n\ge0}\mathcal{M}_{\lambda_{\mathcal{M}}+\frac{n}{l}}=\bigoplus_{n\ge0}\mathcal{M}_{\frac{\lambda_{1}+\cdots+\lambda_{d}}{l}+\frac{d\left(l^{2}-1\right)c}{24l}+\frac{n}{l}},
\]
where
\[\mathcal{M}_{\frac{\lambda_{1}+\cdots+\lambda_{d}}{l}+\frac{d\left(l^{2}-1\right)c}{24l}+\frac{n}{l}}=\sum_{n_1+n_2\cdots +n_d=n}M_{\lambda_{1}+n_{1}}^{1}\otimes M_{\lambda_{2}+n_{2}}^{2}\otimes\cdots\otimes M_{\lambda_{d}+n_{d}}^{d}.\]

\subsubsection{$g$-action on $g^{s}$-twisted modules}

Let $\mathcal{M}=T_{g^{s}}^{M^{1},\cdots,M^{d}}$ be as before. From
the discussion in Section 2.2, $g^s$ acts on $\mathcal{M}$ as follows:
\[
g^{s}\left(w_{1}\otimes w_{2}\otimes\cdots\otimes w_{d}\right)=e^{2\pi in/l}\left(w_{1}\otimes w_{2}\otimes\cdots\otimes w_{d}\right).
\]
for $w_{1}\otimes w_{2}\otimes\cdots\otimes w_{d}\in\mathcal{M}_{\lambda_{\mathcal{M}}+\frac{n}{l}}$
where $w_{i}\in\left(T_{h_{i}^{m}}\left(M^{i}\right)\right)\left(\frac{n_{i}}{l}\right)=M^{i}\left(n_{i}\right)$
with $n_{1}+\cdots+n_{d}=n.$
This suggests us to define
\begin{equation}
g\left(w_{1}\otimes w_{2}\otimes\cdots\otimes w_{d}\right)=e^{\frac{2\pi ixn}{k}}\left(w_{d}\otimes w_{1}\otimes\cdots\otimes w_{d-1}\right)\label{action of g on g^s twisted module}
\end{equation}
by noting that  $T_{g^{s}}^{M_{1},\cdots,M_{d}}\circ g\cong T_{g^{s}}^{M_{d},M_{1},\cdots,M_{d-1}}$
where $x$ is an integer satisfying $sx+ky=d$ for some $y.$

\section{Trace functions and the $S$-Matrix of $V^{\otimes k}$\label{trace functions}}

In this section, we will compute trace functions of twisted $V^{\otimes k}$-modules,
which will be used to compute the $S$-matrix of $V^{\otimes k}$.

\subsection{Trace Functions}

Let $\mathfrak{U}\left(g^{s},g^{r}\right)$ be as given in Lemma \ref{g^s, g^r modules}\emph{.
}In the following, we denote
\[
T_{g^{s}}^{M^{1},\cdots,M^{f};b}=\left(T_{h_{1}^{m}}^{l}\left(M^{1}\right)\right)^{\otimes b}\otimes\cdots\otimes\left(T_{h_{f}^{m}}^{l}\left(M^{f}\right)\right)^{\otimes b}\in\mathfrak{U}\left(g^{s},g^{r}\right)
.\]

\begin{lemma} \label{trace function of g^s-module} Suppose $s,r\in\mathbb{N}$
with $s$ positive. Let $d=\text{gcd}\left(s,k\right)$, $f=\text{gcd}$$\left(d,r\right)$,
$l=\frac{k}{d}$ and $b=\frac{d}{f}$. For $1\le i\le f$, let $\overline{v_{i}}=\sum_{j=1}^{l}v_{i}^{j}$
where $v_{i}$ is a highest weight vector for the Virasoro algebra.
Then
\begin{align*}
 & Z_{T_{g^{s}}^{M^{1},\cdots,M^{f};b}}\left(\overline{v_{1}}\otimes\cdots\otimes\overline{v_{f}},\left(g^{s},g^{r}\right),\tau\right)\\
= & l^{-\left(\text{wt}v_{1}+\cdots+\text{wt}v_{f}\right)+f}e^{-\frac{2\pi ixr}{fl}\left(\lambda_{1}+\cdots+\lambda_{f}-\frac{fc}{24}\right)}Z_{M^{1}}\left(v_{1},\frac{d\tau+rx}{fl}\right)\cdots Z_{M^{f}}\left(v_{f},\frac{d\tau+rx}{fl}\right),
\end{align*}
where $x\in\mathbb{Z}$ satisfies $sx\equiv d$ (mod $k$) and $\lambda_{i}$
is the conformal weight of the irreducible $V$-module $M^{i}$, $1\le i\le f$.

\end{lemma} \begin{proof} For convenience here we denote $\mathcal{M}=T_{g^{s}}^{M^{1},\cdots,M^{f};b}.$
It is clear from (\ref{virasoro of M}) and (\ref{weight of M}) that
\[
L_{\mathcal{M}}\left(0\right)=\sum_{i=1}^{d}\frac{1}{l}L^{i}\left(0\right)+\frac{d\left(l^{2}-1\right)c}{24l}
\]and
\begin{gather*}
\lambda_{\mathcal{M}}=\frac{b\left(\lambda_{1}+\cdots+\lambda_{f}\right)}{l}+\frac{d\left(l^{2}-1\right)c}{24l},
\end{gather*}
where $\lambda_{i}$ is the conformal weight of $M^{i},$ $1\le i\le f$.
Using (\ref{action of g on g^s twisted module}), we have
\begin{alignat*}{1}
 & Z_{\mathcal{M}}\left(\overline{v_{1}}\otimes\cdots\otimes\overline{v_{f}},\left(g^{s},g^{r}\right),\tau\right)\\
= & \text{tr}_{\mathcal{M}}o\left(\overline{v_{1}}\otimes\cdots\otimes\overline{v_{f}}\right)g^{r}q^{L_{\mathcal{M}}\left(0\right)-\frac{kc}{24}}\\
= & \text{tr}_{\mathcal{M}}o\left(\overline{v_{1}}\otimes\cdots\otimes\overline{v_{f}}\right)g^{r}q^{\sum_{i=1}^{d}\frac{1}{l}\left(L^{i}\left(0\right)-\frac{c}{24}\right)}\\
= & \text{\ensuremath{\sum_{n\ge0}}tr}_{\left(T_{h_{1}}^{l}\left(M^{1}\right)\right)_{\frac{\lambda_{1}+n}{l}}^{\otimes b}}o\left(\overline{v_{1}}\right)e^{\frac{2\pi ixnrb}{k}}q^{\sum_{i=1}^{b}\frac{1}{l}\left(L^{i}\left(0\right)-\frac{c}{24}\right)}\\
 & \cdots\text{\ensuremath{\sum_{n\ge0}}tr}_{\left(T_{h_{f}}^{l}\left(M^{f}\right)\right)_{\frac{\lambda_{f}+n}{l}}^{\otimes b}}o\left(\overline{v_{f}}\right)e^{\frac{2\pi ixnrb}{k}}q^{\sum_{i=1}^{b}\frac{1}{l}\left(L^{i}\left(0\right)-\frac{c}{24}\right)}\\
= & \sum_{n\ge0}\text{tr}_{\left(M^{1}\right)_{\lambda_{1}+n}^{\otimes b}}l^{-\text{wt}v_{1}+1}o\left(v_{1}\right)e^{\frac{2\pi ixnr}{fl}}\left(e^{2\pi i\tau}\right)^{\frac{b}{l}\left(\lambda_{1}+n-\frac{c}{24}\right)}\\
 & \cdots\sum_{n\ge0}\text{tr}_{\left(M^{f}\right)_{\lambda_{f}+n}^{\otimes b}}l^{-\text{wt}v_{f}+1}o\left(v_{f}\right)e^{\frac{2\pi ixnr}{fl}}\left(e^{2\pi i\tau}\right)^{\frac{b}{l}\left(\lambda_{f}+n-\frac{c}{24}\right)}\\
= & \sum_{n\ge0}\text{tr}_{\left(M^{1}\right)_{\lambda_{1}+n}^{\otimes b}}l^{-\text{wt}v_{1}+1}o\left(v_{1}\right)\left(e^{\frac{xr}{bf}}e^{\tau}\right)^{\frac{2\pi ib}{l}\left(\lambda_{1}+n-\frac{c}{24}\right)}e^{-\frac{xr}{bf}\cdot\frac{2\pi ib}{l}\left(\lambda_{1}-\frac{c}{24}\right)}\\
 & \cdots\sum_{n\ge0}\text{tr}_{\left(M^{f}\right)_{\lambda_{f}+n}^{\otimes b}}l^{-\text{wt}v_{f}+1}o\left(v_{f}\right)\left(e^{\frac{xr}{bf}}e^{\tau}\right)^{\frac{2\pi ib}{l}\left(\lambda_{f}+n-\frac{c}{24}\right)}e^{-\frac{xr}{bf}\cdot\frac{2\pi ib}{l}\left(\lambda_{f}-\frac{c}{24}\right)}\\
= & l^{-\text{wt}v_{1}+1}e^{-\frac{2\pi ixr}{fl}\left(\lambda_{1}-\frac{c}{24}\right)}Z_{M^{1}}\left(v_{1},\left(\tau+\frac{rx}{bf}\right)\frac{b}{l}\right)\\
 & \cdots l^{-\text{wt}v_{f}+1}e^{-\frac{2\pi ixr}{fl}\left(\lambda_{f}-\frac{c}{24}\right)}Z_{M^{f}}\left(v_{f},\left(\tau+\frac{rx}{bf}\right)\frac{b}{l}\right)\\
= & l^{-\left(\text{wt}v_{1}+\cdots+\text{wt}v_{f}\right)+f}e^{-\frac{2\pi ixr}{fl}\left(\lambda_{1}+\cdots+\lambda_{f}-\frac{fc}{24}\right)}Z_{M^{1}}\left(v_{1},\frac{d\tau+rx}{fl}\right)\cdots Z_{M^{f}}\left(v_{f},\frac{d\tau+rx}{fl}\right).
\end{alignat*}
\end{proof} 

Let $\mathfrak{U}\left(g^{r},g^{s}\right)$ be the equivalence classes
of irreducible $g^{r}$-twisted modules that are $g^{s}$-stable.
Assume gcd$\left(r,k\right)=d_{1}$, then there exist $p,q\in\mathbb{Z}$
such that $rp+kq=d_{1}.$ Also we note that gcd$\left(d_{1},s\right)=\text{gcd}\left(d,r\right)=\text{gcd}\left(s,k,r\right)=f$.
Set $l_{1}=\frac{k}{d_{1}}$ and $a=\frac{d_{1}}{f}$. By similar
arguments as above, we see that any $g^{r}$-twisted module in $\mathfrak{U}\left(g^{r},g^{s}\right)$
can be written in the form
\begin{equation}
T_{g^{r}}^{M^{1},\cdots,M^{f};a}=\left(T_{c_{1}}^{l_{1}}\left(M^{1}\right)\right)^{\otimes a}\otimes\cdots\otimes\left(T_{c_{f}}^{l_{1}}\left(M^{f}\right)\right)^{\otimes a}\label{general module g^r,g^s},
\end{equation}
where $M^{1},\cdots,M^{d}\in\text{Irr}\left(V\right)$ and each $c_{i}$
is an $l_{1}$-cycle, $1\le i\le f$. By Lemma \ref{trace function of g^s-module}
we have the following trace function:
\begin{align}
 & Z_{T_{g^{r}}^{M^{1},\cdots,M^{f};a}}\left(\overline{v_{1}}\otimes\cdots\otimes\overline{v_{f}},\left(g^{r},g^{s}\right),\tau\right)\nonumber \\
= & l_{1}^{-\left(\text{wt}v_{1}+\cdots+\text{wt}v_{f}\right)+f}e^{-\frac{2\pi ips}{fl_{1}}\left(\lambda_{1}+\cdots+\lambda_{f}-\frac{fc}{24}\right)}Z_{M^{1}}\left(v_{1},\frac{d_{1}\tau+sp}{fl_{1}}\right)\cdots Z_{M^{f}}\left(v_{f},\frac{d_{1}\tau+sp}{fl_{1}}\right),\label{eq:g^r,s}
\end{align}
where $p\in\mathbb{Z}$ satisfies $rp\equiv d_{1}$ $\left(\text{mod}\ k\right)$,
$\overline{v_{i}}=\sum_{j=1}^{l}v_{i}^{j}$ and $\lambda_{i}$ is
the weight of $M^{i}$, $1\le i\le f.$

Now we consider trace functions of modules in $\mathcal{\mathfrak{U}}\left(1,g^{s}\right).$
Let $i_{1},\cdots,i_{k}\in\left\{ 0,\cdots,p\right\} $. Set $W^{i_{1},\cdots,i_{k}}=W^{i_{1}}\otimes\cdots\otimes W^{i_{k}}$.
Then $W^{i_{1},\cdots,i_{k}}$ is an irreducible $V^{\otimes k}$-module.
Recall from (\ref{cycle decompositon of g^s}) that $g^{s}$ can be
written as a product of $d$ disjoint $l$-cycles: $g^{s}=h_{1}^{m}\cdots h_{d}^{m}$
where $d=\text{gcd}\left(s,k\right)$, $s=dm$ and $k=dl$. Therefore,
$W^{i_{1},\cdots,i_{k}}\circ g^{s}\cong W^{i_{1},\cdots,i_{k}}$ if
and only if $W^{i_a}=W^{i_{a}+jd}$ for $i_a=1, \cdots, d$ and $j=1,\cdots, l-1.$ Therefore, any $V^{\otimes k}$-module
in $\mathcal{\mathfrak{U}}\left(1,g^{s}\right)$ can be written in
the form
\begin{equation}
\left(W^{i_{1}}\right)^{\otimes l}\otimes\cdots\otimes\left(W^{i_{d}}\right)^{\otimes l},\label{1,g^s}
\end{equation}
where $i_{1},\cdots,i_{d}\in\left\{ 0,1,\cdots,p\right\} $.

We will need to use the following results that were proved in \cite{DRX21}:

\begin{lemma} \label{trace function of tensor products of V-modules}Let
$M$ be an irreducible $V$-module and $h$ be an $n$-cycle, where
$n$ is a positive integer. Then

(1) $Z_{M^{\otimes n}}\left(\boldsymbol{1},\left(1,h\right),\tau\right)=\chi_{M}\left(n\tau\right)$;

(2) $Z_{M^{\otimes n}}\left(\bar{v},\left(1,h\right),\tau\right)=nZ_{M}\left(v,n\tau\right)$
where $\bar{v}=\sum_{j=1}^{n}v^{j}$ with $v\in V$.

\end{lemma}

The following trace function will be used later to find entries involving
untwisted and twisted $V^{\otimes k}$-modules in the $S$-matrix.

\begin{lemma} \label{trace function of (1,g^s)-modules} Suppose
that $s\in\mathbb{N}$, $d=\text{gcd}\left(s,k\right)$, $m=\frac{s}{d}$
and $l=\frac{k}{d}.$ For $1\le i\le d$, let $\overline{v_{i}}=\sum_{j=1}^{l}v_{i}^{j}$
where each $v_{i}$ is a highest weight vector for the Virasoro algebra.
Then
\[
Z_{\left(W^{i_{1}}\right)^{\otimes l}\otimes\cdots\otimes\left(W^{i_{d}}\right)^{\otimes l}}\left(\overline{v_{1}}\otimes\cdots\otimes\overline{v_{d}},\left(1,g^{s}\right),\tau\right)=l^{d}Z_{W^{i_{1}}}\left(v_{1},l\tau\right)\cdots Z_{W^{i_{d}}}\left(v_{d},l\tau\right).\]
\end{lemma}
\begin{proof} A straightforward calculation gives
\begin{alignat*}{1}
 & Z_{\left(W^{i_{1}}\right)^{\otimes l}\otimes\cdots\otimes\left(W^{i_{d}}\right)^{\otimes l}}\left(\overline{v_{1}}\otimes\cdots\otimes\overline{v_{d}},\left(1,g^{s}\right),\tau\right)\\
= & Z_{\left(W^{i_{1}}\right)^{\otimes l}\otimes\cdots\otimes\left(W^{i_{d}}\right)^{\otimes l}}\left(\overline{v_{1}}\otimes\cdots\otimes\overline{v_{d}},\left(1,h_{1}^{m}\cdots h_{d}^{m}\right),\tau\right)\\
= & Z_{\left(W^{i_{1}}\right)^{\otimes l}}\left(\overline{v_{1}},\left(1,h_{1}^{m}\right),\tau\right)\cdots Z_{\left(W^{i_{d}}\right)^{\otimes l}}\left(\overline{v_{d}},\left(1,h_{d}^{m}\right),\tau\right)\\
= & lZ_{W^{i_{1}}}\left(v_{1},l\tau\right)\cdots lZ_{W^{i_{d}}}\left(v_{d},l\tau\right)\\
= & l^{d}Z_{W^{i_{1}}}\left(v_{1},l\tau\right)\cdots Z_{W^{i_{d}}}\left(v_{d},l\tau\right),
\end{alignat*}
where we use (\ref{cycle decompositon of g^s}) and Lemma \ref{trace function of tensor products of V-modules}.
\end{proof}

\subsection{$S$-matrix for $V^{\otimes k}$ }

In this subsection, we will apply results in Subsection 4.1 to obtain
the $S$-matrix of $V^{\otimes k}.$

Recall that $\text{Irr}\left(V\right)=\left\{ W^{0},\cdots,W^{p}\right\} $
is the set of all inequivalent irreducible $V$-modules. Let $a,b,l,l_{1},f$
be as before and $i_{1},\cdots,i_{f},j_{1},\cdots,j_{f}\in\left\{ 0,1,\cdots,p\right\} $.
For convenience, in the following we will denote
\begin{gather*}
T_{g^{r}}^{i_{1},\cdots,i_{f};a}=\left(T_{c_{1}}^{l_{1}}\left(W^{i_{1}}\right)\right)^{\otimes a}\otimes\cdots\otimes\left(T_{c_{f}}^{l_{1}}\left(W^{i_{f}}\right)\right)^{\otimes a}\in\mathfrak{U}\left(g^{r},g^{s}\right),\\
T_{g^{s}}^{j_{1},\cdots,j_{f};b}=\left(T_{h_{1}^{m}}^{l}\left(W^{j_{1}}\right)\right)^{\otimes b}\otimes\cdots\otimes\left(T_{h_{f}^{m}}^{l}\left(W^{j_{f}}\right)\right)^{\otimes b}\in\mathfrak{U}\left(g^{s},g^{r}\right),
\end{gather*}
where $\mathfrak{U}\left(g^{r},g^{s}\right)$ and $\mathfrak{U}\left(g^{s},g^{r}\right)$
are given in (\ref{general module g^r,g^s}) and (\ref{general module g^s,g^r})
respectively. By abusing the notations, we will also denote the conformal
weight of $W^{i}$ by $\lambda_{i}$ for  $0\le i\le p$.

\begin{lemma} \label{S-matrix: twisted and twisted} Suppose $s,r\in\mathbb{N}$
with $s$ positive. Let $d=\text{gcd}\left(s,k\right)$, $d_{1}=\text{gcd}\left(r,k\right)$,
$f=\text{gcd}$$\left(d,r\right)$. Then there exist  $x,y,p,q\in\mathbb{Z}$ such that
 $sx+ky=d$ and $pr+kq=d_{1}.$ Set $l=\frac{k}{d}$, $l_{1}=\frac{k}{d_{1}},$ $b=\frac{d}{f}$
and $a=\frac{d_{1}}{f}$. Let $T_{g^{r}}^{i_{1},\cdots,i_{f};a}\in\mathfrak{U}\left(g^{r},g^{s}\right)$
and $T_{g^{s}}^{j_{1},\cdots,j_{f};b}\in\mathfrak{U}\left(g^{s},g^{r}\right)$
as above. Then
\begin{alignat*}{1}
S_{T_{g^{r}}^{i_{1},\cdots,i_{f};a},T_{g^{s}}^{j_{1},\cdots,j_{f};b}} & =\left(\frac{l_{1}}{l}\right)^{f}e^{-\frac{2\pi ips}{fl_{1}}\left(\lambda_{i_{1}}+\cdots+\lambda_{i_{f}}-\frac{fc}{24}\right)}e^{-\frac{2\pi ixr}{fl}\left(\lambda_{j_{1}}+\cdots+\lambda_{j_{f}}-\frac{fc}{24}\right)}\\
{\color{red}} & \cdot\sum_{n_{1}=0}^{p}S_{W^{i_{1}},W^{n_{1}}}A_{n_{1},j_{1}}^{r,s}\cdots\sum_{n_{f}=0}^{p}S_{W^{i_{f}},W^{n_{f}}}A_{n_{f},j_{f}}^{r,s}
\end{alignat*}
where $A_{n_{t},j_{t}}^{r,s}$is the entry of $\rho\left(A^{r,s}\right)$
defined in Theorem \ref{modular invariance thm} with $A^{r,s}=\left[\begin{array}{cc}
\frac{l_{1}}{b} & \frac{rx}{d_{1}}\\
\frac{-sp}{d} & \frac{dq+yd_{1}-yqk}{f}
\end{array}\right]$, and $\lambda_{i_{t}}$ and $\lambda_{j_{t}}$ are the conformal
weights of $W^{i_{t}}$ and $W^{j_{t}}$, respectively, for $1\le t\le f$.

\end{lemma} \begin{proof} First we prove that $A^{r,s}\in SL_{2}\left(\mathbb{Z}\right).$
Indeed, since $k=dl=d_{1}l_{1}$ with $d=bf$ and $d_{1}=af$, we
obtain $bl=al_{1}.$ Since gcd$\left(a,b\right)=1,$ we see that  $b\mid l_{1}$ and  $\frac{l_{1}}{b}\in\mathbb{Z}$.
The entries $\frac{rx}{d_{1}}$ and $\frac{-sp}{d}$ are obviously
in $\mathbb{Z}$ since $d_{1}\mid r$ and $d\mid s$. Since $f=\text{gcd}\left(d,r\right)=\text{gcd}\left(d_{1},s\right)$
where $d=\text{gcd}\left(s,k\right)$ and $d_{1}=\text{gcd}\left(r,k\right)$,
it is clear that $\frac{dq+yd_{1}-yqk}{f}\in\mathbb{Z}.$ It is easy
to check that $\det(A^{r,s})=1.$

Let $v_{1},\cdots,v_{f}\in V$ be highest weight vectors for the Virasoro
algebra and $\overline{v_{i}}=\sum_{j=1}^{l}v_{i}^{j},1\le i\le f$.
By (\ref{eq:g^r,s}) we have
\begin{alignat*}{1}
 & Z_{T_{g^{r}}^{i_{1},\cdots,i_{f};a}}\left(\overline{v_{1}}\otimes\cdots\otimes\overline{v_{f}},\left(g^{r},g^{s}\right),-1/\tau\right)\\
= & l_{1}^{-\left(\text{wt}v_{1}+\cdots+\text{wt}v_{f}\right)+f}e^{-\frac{2\pi ips}{fl_{1}}\left(\lambda_{i_{1}}+\cdots+\lambda_{i_{f}}-\frac{fc}{24}\right)}Z_{W^{i_{1}}}\left(v_{1},\frac{-d_{1}/\tau+sp}{fl_{1}}\right)\cdots Z_{W^{i_{f}}}\left(v_{f},\frac{-d_{1}/\tau+sp}{fl_{1}}\right)\\
= & l_{1}^{-\left(\text{wt}v_{1}+\cdots+\text{wt}v_{f}\right)+f}e^{-\frac{2\pi ips}{fl_{1}}\left(\lambda_{i_{1}}+\cdots+\lambda_{i_{f}}-\frac{fc}{24}\right)}Z_{W^{i_{1}}}\left(v_{1},-\frac{1}{\frac{fl_{1}\tau}{d_{1}-sp\tau}}\right)\cdots Z_{W^{i_{f}}}\left(v_{f},-\frac{1}{\frac{fl_{1}\tau}{d_{1}-sp\tau}}\right)\\
= & l_{1}^{-\left(\text{wt}v_{1}+\cdots+\text{wt}v_{f}\right)+f}e^{-\frac{2\pi ips}{fl_{1}}\left(\lambda_{i_{1}}+\cdots+\lambda_{i_{f}}-\frac{fc}{24}\right)}Z_{W^{i_{1}}}\left(v_{1},-\frac{1}{\bar{\tau}}\right)\cdots Z_{W^{i_{f}}}\left(v_{f},-\frac{1}{\bar{\tau}}\right)\\
= & l_{1}^{-\left(\text{wt}v_{1}+\cdots+\text{wt}v_{f}\right)+f}e^{-\frac{2\pi ips}{fl_{1}}\left(\lambda_{i_{1}}+\cdots+\lambda_{i_{f}}-\frac{fc}{24}\right)}\left(\bar{\tau}^{\text{wt}v_{1}}\right)\sum_{i=0}^{p}S_{W^{i_{1}},W^{i}}Z_{W^{i}}\left(v_{1},\bar{\tau}\right)\cdots\left(\bar{\tau}^{\text{wt}v_{f}}\right)\sum_{i=0}^{p}S_{W^{i_{f}},W^{i}}Z_{W^{i}}\left(v_{f},\bar{\tau}\right)
\end{alignat*}
where $\bar{\tau}=\frac{fl_{1}\tau}{d_{1}-sp\tau}=A^{r,s}\frac{d\tau-rx}{fl}$
with $A^{r,s}=\left[\begin{array}{cc}
\frac{fl_{1}}{d} & \frac{rxl_{1}}{k}\\
\frac{-sp}{d} & \frac{dd_{1}-sxrp}{kf}
\end{array}\right]=\left[\begin{array}{cc}
\frac{l_{1}}{b} & \frac{rx}{d_{1}}\\
\frac{-sp}{d} & \frac{dq+yd_{1}-yqk}{f}
\end{array}\right]\in SL_{2}\left(\mathbb{Z}\right).$

Note that for $v_{t}\in V$, $1\le t\le f$, we have
\[
Z_{W^{i}}\left(v_{t},\bar{\tau}\right)=\left(\frac{d_{1}-sp\tau}{fl}\right)^{\text{wt}v_{t}}\sum_{n=0}^{p}A_{i,n}^{r,s}Z_{W^{n}}\left(v_{t},\frac{d\tau-rx}{fl}\right),0\le i\le p.
\]
Thus we have
\begin{alignat}{1}
 & Z_{T_{g^{r}}^{i_{1},\cdots,i_{f};a}}\left(\overline{v_{1}}\otimes\cdots\otimes\overline{v_{f}},\left(g^{r},g^{s}\right),-1/\tau\right)\nonumber \\
= & l_{1}^{-\left(\text{wt}v_{1}+\cdots+\text{wt}v_{f}\right)+f}e^{-\frac{2\pi ips}{fl_{1}}\left(\lambda_{i_{1}}+\cdots+\lambda_{i_{f}}-\frac{fc}{24}\right)}\nonumber \\
 & \cdot\left(\bar{\tau}^{\text{wt}v_{1}}\right)\sum_{i=0}^{p}S_{W^{i_{1}},W^{i}}\left(\frac{d_{1}-sp\tau}{fl}\right)^{\text{wt}v_{1}}\sum_{n=0}^{p}A_{i,n}^{r,s}Z_{W^{n}}\left(v_{1},\frac{d\tau-rx}{fl}\right)\nonumber \\
 & \cdots\left(\bar{\tau}^{\text{wt}v_{f}}\right)\sum_{i=0}^{p}S_{W^{i_{f}},W^{i}}\left(\frac{d_{1}-sp\tau}{fl}\right)^{\text{wt}v_{f}}\sum_{n=0}^{p}A_{i,n}^{r,s}Z_{W^{n}}\left(v_{f},\frac{d\tau-rx}{fl}\right)\nonumber \\
= & l_{1}^{-\left(\text{wt}v_{1}+\cdots+\text{wt}v_{f}\right)+f}e^{-\frac{2\pi ips}{fl_{1}}\left(\lambda_{i_{1}}+\cdots+\lambda_{i_{f}}-\frac{fc}{24}\right)}\nonumber \\
 & \cdot\sum_{i=0}^{p}S_{W^{i_{1}},W^{i}}\left(\frac{l_{1}\tau}{l}\right)^{\text{wt}v_{1}}\sum_{n=0}^{p}A_{i,n}^{r,s}Z_{W^{n}}\left(v_{1},\frac{d\tau-rx}{fl}\right)\nonumber \\
 & \cdots\sum_{i=0}^{p}S_{W^{i_{f}},W^{i}}\left(\frac{l_{1}\tau}{l}\right)^{\text{wt}v_{f}}\sum_{n=0}^{p}A_{i,n}^{r,s}Z_{W^{n}}\left(v_{f},\frac{d\tau-rx}{fl}\right)\nonumber \\
= & l_{1}^{f}e^{-\frac{2\pi ips}{fl_{1}}\left(\lambda_{i_{1}}+\cdots+\lambda_{i_{f}}-\frac{fc}{24}\right)}l^{-\left(\text{wt}v_{1}+\cdots+\text{wt}v_{f}\right)}\tau^{\text{wt}v_{1}+\cdots+\text{wt}v_{f}}\nonumber \\
 & \cdot\sum_{j_{1},n_{1}=0}^{p}S_{W^{i_{1}},W^{j_{1}}}A_{j_{1},n_{1}}^{r,s}Z_{W^{n_{1}}}\left(v_{1},\frac{d\tau-rx}{fl}\right)\cdots\sum_{j_{f},n_{f}=0}^{p}S_{W^{i_{f}},W^{j_{f}}}A_{j_{f},n_{f}}^{r,s}Z_{W^{n_{f}}}\left(v_{f},\frac{d\tau-rx}{fl}\right)\nonumber \\
= & l_{1}^{f}e^{-\frac{2\pi ips}{fl_{1}}\left(\lambda_{i_{1}}+\cdots+\lambda_{i_{f}}-\frac{fc}{24}\right)}l^{-\left(\text{wt}v_{1}+\cdots+\text{wt}v_{f}\right)}\tau^{\text{wt}v_{1}+\cdots+\text{wt}v_{f}}\nonumber \\
 & \cdot\sum_{j_{1},n_{1}=0}^{p}S_{W^{i_{1}},W^{n_{1}}}A_{n_{1},j_{1}}^{r,s}Z_{W^{j_{1}}}\left(v_{1},\frac{d\tau-rx}{fl}\right)\cdots\sum_{j_{f},n_{f}=0}^{p}S_{W^{i_{f}},W^{n_{f}}}A_{n_{f},j_{f}}^{r,s}Z_{W^{j_{f}}}\left(v_{f},\frac{d\tau-rx}{fl}\right).\label{compare 1}
\end{alignat}

On the other hand,
\begin{alignat}{1}
 & Z_{T_{g^{r}}^{i_{1},\cdots,i_{f};a}}\left(\overline{v_{1}}\otimes\cdots\otimes\overline{v_{f}},\left(g^{r},g^{s}\right),-1/\tau\right)\nonumber \\
= & \tau^{\text{wt}v_{1}+\cdots+\text{wt}v_{f}}\cdot\sum_{\left\{ j_{1},\cdots,j_{f}\right\} \subset\left\{ 0,1,\cdots,p\right\} }S_{T_{g^{r}}^{i_{1},\cdots,i_{f};a},T_{g^{s}}^{j_{1},\cdots,j_{f};b}}Z_{T_{g^{s}}^{j_{1},\cdots,j_{f};b}}\left(\overline{v_{1}}\otimes\cdots\otimes\overline{v_{f}},\left(g^{s},g^{-r}\right),\tau\right)\nonumber \\
= & \tau^{\text{wt}v_{1}+\cdots+\text{wt}v_{f}}\cdot\sum_{\left\{ j_{1},\cdots,j_{f}\right\} \subset\left\{ 0,1,\cdots,p\right\} }S_{T_{g^{r}}^{i_{1},\cdots,i_{f};a},T_{g^{s}}^{j_{1},\cdots,j_{f};b}}l^{-\left(\text{wt}v_{1}+\cdots+\text{wt}v_{f}\right)+f}e^{\frac{2\pi ixr}{fl}\left(\lambda_{j_{1}}+\cdots+\lambda_{j_{f}}-\frac{fc}{24}\right)}\nonumber \\
{\color{red}} & \cdot Z_{W^{j_{1}}}\left(v_{1},\frac{d\tau-rx}{fl}\right)\cdots Z_{W^{j_{f}}}\left(v_{f},\frac{d\tau-rx}{fl}\right),\label{compare 2}
\end{alignat}
where we use Lemma \ref{trace function of g^s-module}. Comparing
(\ref{compare 1}) and (\ref{compare 2}), we obtain
\begin{alignat*}{1}
S_{T_{g^{r}}^{i_{1},\cdots,i_{f};a},T_{g^{s}}^{j_{1},\cdots,j_{f};b}} & =\left(\frac{l_{1}}{l}\right)^{f}e^{-\frac{2\pi ips}{fl_{1}}\left(\lambda_{i_{1}}+\cdots+\lambda_{i_{f}}-\frac{fc}{24}\right)}e^{-\frac{2\pi ixr}{fl}\left(\lambda_{j_{1}}+\cdots+\lambda_{j_{f}}-\frac{fc}{24}\right)}\\
{\color{red}} & \cdot\sum_{n_{1}=0}^{p}S_{W^{i_{1}},W^{n_{1}}}A_{n_{1},j_{1}}^{r,s}\cdots\sum_{n_{f}=0}^{p}S_{W^{i_{f}},W^{n_{f}}}A_{n_{f},j_{f}}^{r,s}.
\end{alignat*}
\end{proof}

Recall from (\ref{1,g^s}) and (\ref{g^s-module general form}) that any $V^{\otimes k}$-module
in $\mathfrak{U}(1,g^{s})$ has the form $\left(W^{i_{1}}\right)^{\otimes l}\otimes\cdots\otimes\left(W^{i_{d}}\right)^{\otimes l}$ with $i_{1},\cdots,i_{d}\in\left\{ 0,1,\cdots,p\right\},$ and  that any $g^{s}$-twisted module has the form
\[
T_{g^{s}}^{W^{j_{1}},\cdots,W^{j_{d}}}=T_{h_{1}^{m}}^{l}\left(W^{j_{1}}\right)\otimes\cdots\otimes T_{h_{d}^{m}}^{l}\left(W^{j_{d}}\right)
\]
with $j_{1},\cdots,j_{d}\in\left\{ 0,1,\cdots,p\right\}.$ Now we
find entries involving untwisted and twisted $V^{\otimes k}$-modules
in the $S$-matrix.

\begin{lemma} \label{S-matrix, nontwisted and twisted} Suppose $s\in\mathbb{N},$
$d=\text{gcd}\left(s,k\right)$ and $l=\frac{k}{d}.$ Let $\left(W^{i_{1}}\right)^{\otimes l}\otimes\cdots\otimes\left(W^{i_{d}}\right)^{\otimes l}\in\mathcal{\mathfrak{U}}\left(1,g^{s}\right)$
and $T_{g^{s}}^{W^{j_{1}},\cdots,W^{j_{d}}}$ be a $g^{s}$-twisted
$V^{\otimes k}$-module as above. Then
\[
S_{\left(W^{i_{1}}\right)^{\otimes l}\otimes\cdots\otimes\left(W^{i_{d}}\right)^{\otimes l},T_{g^{s}}^{W^{j_{1}},\cdots,W^{j_{d}}}}=S_{W^{i_{1}},W^{j_{1}}}\cdots S_{W^{i_{d}},W^{j_{d}}}.
\]

\end{lemma} \begin{proof} For $1\le i\le d$, let $\overline{v_{i}}=\sum_{j=1}^{l}v_{i}^{j}$
where each $v_{i}$ is a highest weight vector for the Virasoro algebra.
By Lemma \ref{trace function of (1,g^s)-modules}, we have
\begin{alignat}{1}
 & Z_{\left(W^{i_{1}}\right)^{\otimes l}\otimes\cdots\otimes\left(W^{i_{d}}\right)^{\otimes l}}\left(\overline{v_{1}}\otimes\cdots\otimes\overline{v_{d}},\left(1,g^{s}\right),-\frac{1}{\tau}\right)\nonumber \\
= & l^{d}Z_{W^{i_{1}}}\left(v_{1},-\frac{l}{\tau}\right)\cdots Z_{W^{i_{d}}}\left(v_{d},-\frac{l}{\tau}\right)\nonumber \\
= & l^{d}Z_{W^{i_{1}}}\left(v_{1},\frac{-1}{\frac{\tau}{l}}\right)\cdots Z_{W^{i_{d}}}\left(v_{d},\frac{-1}{\frac{\tau}{l}}\right)\nonumber \\
= & l^{d}\left(\frac{\tau}{l}\right)^{\text{wt}v_{1}}\sum_{j_{1}=0}^{p}S_{W^{i_{1}},W^{j_{1}}}Z_{W^{j_{1}}}\left(v_{1},\frac{\tau}{l}\right)\cdots\left(\frac{\tau}{l}\right)^{\text{wt}v_{d}}\sum_{j_{d}=0}^{p}S_{W^{i_{d}},W^{j_{d}}}Z_{W^{j_{d}}}\left(v_{d},\frac{\tau}{l}\right)\nonumber \\
= & \tau^{\text{wt}v_{1}+\cdots+\text{wt}v_{d}}l^{-\left(\text{wt}v_{1}+\cdots+\text{wt}v_{d}\right)+d}\sum_{j_{1}=0}^{p}S_{W^{i_{1}},W^{j_{1}}}Z_{W^{j_{1}}}\left(v_{1},\frac{\tau}{l}\right)\cdots\sum_{j_{d}=0}^{p}S_{W^{i_{d}},W^{j_{d}}}Z_{W^{j_{d}}}\left(v_{d},\frac{\tau}{l}\right).\label{(1,g^s)-1}
\end{alignat}
By Lemma \ref{trace function of g^s-module} with $r=k$ and $f=d$,
we have
\begin{align*}
 & Z_{T_{g^{s}}^{W^{j_{1}},\cdots,W^{j_{d}}}}\left(\overline{v_{1}}\otimes\cdots\otimes\overline{v_{d}},\left(g^{s},1\right),\tau\right)\\
= & l^{-\left(\text{wt}v_{1}+\cdots+\text{wt}v_{d}\right)+d}Z_{W^{j_{1}}}\left(v_{1},\frac{\tau-lx}{l}\right)\cdots Z_{W^{j_{d}}}\left(v_{d},\frac{\tau-lx}{l}\right)\\
= & l^{-\left(\text{wt}v_{1}+\cdots+\text{wt}v_{d}\right)+d}Z_{W^{j_{1}}}\left(v_{1},\frac{\tau}{l}\right)\cdots Z_{W^{j_{d}}}\left(v_{d},\frac{\tau}{l}\right),
\end{align*}
where $x\in\mathbb{Z}$ satisfies $sx\equiv d$ (mod $k$).

On the other hand, we have
\begin{alignat}{1}
 & Z_{\left(W^{i_{1}}\right)^{\otimes l}\otimes\cdots\otimes\left(W^{i_{d}}\right)^{\otimes l}}\left(\overline{v_{1}}\otimes\cdots\otimes\overline{v_{d}},\left(1,g^{s}\right),-\frac{1}{\tau}\right)\nonumber \\
= & \tau^{\text{wt}v_{1}+\cdots+\text{wt}v_{d}}\sum_{\left\{ j_{1},\cdots,j_{p}\right\} \subset\left\{ 0,1,\cdots,p\right\} }S_{\left(W^{i_{1}}\right)^{\otimes l}\otimes\cdots\otimes\left(W^{i_{d}}\right)^{\otimes l},T_{g^{s}}^{W^{j_{1},}\cdots,W^{j_{d}}}}Z_{T_{g^{s}}^{W^{j_{1}},\cdots,W^{j_{d}}}}\left(\overline{v_{1}}\otimes\cdots\overline{v_{d}},\left(g^{s},1\right),\tau\right)\nonumber \\
= & \tau^{\text{wt}v_{1}+\cdots+\text{wt}v_{d}}\sum_{\left\{ j_{1},\cdots,j_{p}\right\} \subset\left\{ 0,1,\cdots,p\right\} }S_{\left(W^{i_{1}}\right)^{\otimes l}\otimes\cdots\otimes\left(W^{i_{d}}\right)^{\otimes l},T_{g^{s}}^{W^{j_{1},}\cdots,W^{j_{d}}}}\nonumber \\
 & \cdot l^{-\left(\text{wt}v_{1}+\cdots+\text{wt}v_{d}\right)+d}Z_{W^{j_{1}}}\left(v_{1},\frac{\tau}{l}\right)\cdots Z_{W^{j_{d}}}\left(v_{d},\frac{\tau}{l}\right).\label{(1,g^s)-2}
\end{alignat}

Comparing the right sides of (\ref{(1,g^s)-1}) and (\ref{(1,g^s)-2}),
we obtain
\[
S_{\left(W^{i_{1}}\right)^{\otimes l}\otimes\cdots\otimes\left(W^{i_{d}}\right)^{\otimes l},T_{g^{s}}^{W^{j_{1}},\cdots,W^{j_{d}}}}=S_{W^{i_{1}},W^{j_{1}}}\cdots S_{W^{i_{d}},W^{j_{d}}}.
\]
\end{proof}

\section{$S$-matrix for $\left(V^{\otimes k}\right)^{\left\langle g\right\rangle }$ }

In this section, we first give a complete list of irreducible modules
for the cyclic permutation orbifold $\left(V^{\otimes k}\right)^{\left\langle g\right\rangle }$.
 Then we will give a precise formula for the  $S$-matrix
of $\left(V^{\otimes k}\right)^{\left\langle g\right\rangle }$ by using representations of $SL_2(\mathbb Z)$ of twisted conformal blocks given in Section \ref{trace functions}.

\subsection{Irreducible modules for $\left(V^{\otimes k}\right)^{\left\langle g\right\rangle }$}

Let $U$ be a vertex operator algebra and $G\le\text{Aut\ensuremath{\left(U\right)}.}$
First we give the irreducible $U^{G}$-modules appearing in an irreducible
$h$-twisted $U$-module for some $h\in G$ \cite{DRX17}. Recall
from Section \ref{modular invariance} that $G$ acts on set $\mathcal{S}=\cup_{h\in G}\mathfrak{U}\left(h\right)$
and $M\circ h$ and $\text{\ensuremath{M}}$ are isomorphic $U^{G}$-modules
for any $h\in G$ and $M\in\mathcal{S}.$ It is obvious that the cardinality
of the $G$-orbit $\left|M\circ G\right|$ of $M$ is equal to $\left[G:G_{M}\right]$
where $G_{M}=\left\{ h\in G\mid M\circ h\cong M\right\} $.

Let $\mathcal{S}=\cup_{j\in J}\mathcal{O}_{j}$ be the decomposition
of $\mathcal{S}$ into a disjoint union of orbits. Let $M^{j}$ for
$j\in J$ be the orbit representatives of $\mathcal{S}$ and $\mathcal{O}_{j}=\left\{ M^{j}\circ h\mid h\in G\right\} $
be the orbit of $M^{j}$ under $G.$ Let $\Lambda_{G_{M}}$ be the
set of all irreducible characters $\lambda$ of $\mathbb{C}^{\alpha_{\mathcal{M}}}\left[G_{M}\right].$
Denote the corresponding simple module by $W_{\lambda}$. Let $M^{\lambda}$
be the sum of simple $\mathbb{C}^{\alpha_{\mathcal{M}}}\left[G_{M}\right]$-submodules
of $M$ isomorphic to $W_{\lambda}.$ Then\[M=\oplus_{\lambda\in\Lambda_{G_{M}}}M^{\lambda}=\oplus_{\lambda\in\Lambda_{G_{M}}}W_{\lambda}\otimes M_{\lambda},
\]where the multiplicity space $M_{\lambda}$ of $W_{\lambda}$ in $M$
is a $U^{G}$-module. The following result is given in \cite{DRX17}
(see also Theorem 4.2 in \cite{DRX21}):

\begin{proposition} \label{abstract complete list}

Suppose that $U$ is a regular vertex operator algebra of CFT type,
$G$ is solvable and the weight of any irreducible twisted $U$-module
$M$ is positive except $U$ itself. Then $\left\{ M_{\lambda}^{j}\mid j\in J,\lambda\in\Lambda_{G_{M^{j}}}\right\} $
gives a complete list of inequivalent irreducible $U^{G}$-modules
appearing in the irreducible twisted $U$-modules. \end{proposition}

If $U=V^{\otimes k}$, $G=\left\langle g\right\rangle $ with $g=\left(1,2,\cdots,k\right)$.
Then every irreducible $\left(V^{\otimes k}\right)^{\left\langle g\right\rangle }$-module
appears in an irreducible $g^{i}$-twisted $V^{\otimes k}$-module
for some $i=0,1,\cdots,k-1.$

Let $j_{1},\cdots,j_{k}\in\left[0,p\right]$ where $\left[0,p\right]=\left\{ 0,1,\cdots,p\right\} $.
As before, set $W^{j_{1},\cdots,j_{k}}=W^{j_{1}}\otimes\cdots\otimes W^{j_{k}}$.
Then $W^{j_{1},\cdots,j_{k}}$ is an irreducible $\left(V^{\otimes k}\right)^{\left\langle g\right\rangle }$-module
if the cardinality of $\left\{ j_{1},\cdots,j_{k}\right\} $ is greater
than 1, and $W^{\sigma\left(j_{1}\right),\cdots,\sigma\left(j_{k}\right)}$
is isomorphic to $W^{j_{1},\cdots,j_{k}}$ for any $\sigma\in S_{k}$.
For any $j\in\left[0,p\right]$,   $W^{j,\cdots,j}$ is a direct sum
of irreducible $\left(V^{\otimes k}\right)^{\left\langle g\right\rangle }$-modules
$\left(W^{j,\cdots,j}\right)^{n}$ for $n=0,\cdots,k-1$ where\[
\left(W^{j,\cdots,j}\right)^{n}=\left\{ w\in W^{j,\cdots,j}\mid gw=e^{\frac{-2\pi in}{k}}w\right\} =\left\{ \sum_{s=0}^{k-1}e^{\frac{2\pi isn}{k}}g^{s}w\mid w\in W^{j,\cdots,j}\right\}
\]and $g$ acts on $W^{j,\cdots,j}$ in an obvious way.

Recall from (\ref{g^s-module general form}) that any $g^{s}$-twisted
$V^{\otimes k}$-module can be written as
\[
T_{g^{s}}^{j_{1},\cdots,j_{d}}=T_{h_{1}^{m}}^{l}\left(W^{j_{1}}\right)\otimes\cdots\otimes T_{h_{d}^{m}}^{l}\left(W^{j_{d}}\right),
\]
where $d=\text{gcd}\left(s,k\right),$ $k=dl$ and $j_{1},\cdots,j_{d}\in\left[0,p\right]$$.$
Using (\ref{g^d stabilizer}), we obtain $G_{T_{g^{s}}^{j_{1},\cdots,j_{d}}}=\left\langle g^{d}\right\rangle $
if $|\{j_{1},\cdots,j_{d}\}|>1$ and $G_{T_{g^{s}}^{j_{1},\cdots,j_{d}}}=\left\langle g\right\rangle $
otherwise.

If $|\{j_{1},\cdots,j_{d}\}|>1$, we have $|G_{T_{g^{s}}^{j_{1},\cdots,j_{d}}}|=o\left(g^{d}\right)=l$
and the cardinality of the orbit of $T_{g^{s}}^{j_{1},\cdots,j_{d}}$
is
\begin{equation}
\left|\mathcal{O}_{T_{g^{s}}^{j_{1},\cdots,j_{d}}}\right|=\left[G:G_{T_{g^{s}}^{j_{1},\cdots,j_{d}}}\right]=d.\label{car of orbit of g^s-twisted module}
\end{equation}
In this case,
\begin{equation}
T_{g^{s}}^{j_{1},\cdots,j_{d}}=\oplus_{t=0}^{l-1}\left(T_{g^{s}}^{j_{1},\cdots,j_{d}}\right)^{t},
\label{g^s-decomposition}\end{equation}
which is a direct sum of irreducible $\left(V^{\otimes k}\right)^{\left\langle g\right\rangle }$-modules
such that $g^{d}$ acts on $\left(T_{g^{s}}^{j_{1},\cdots,j_{d}}\right)^{t}$
as $e^{\frac{-2\pi it}{l}}$ where $0\le t<l$. Notice that $\Lambda_{G_{T_{g^{s}}^{j_{1},\cdots,j_{d}}}}=\left\{ \lambda_{t}\mid t=0,\cdots,l-1\right\} $
where $\lambda_{t}\left(g^{d}\right)=e^{-\frac{2\pi it}{l}}$. If
$d=1$, then it reduces to the case given in \cite{DRX21}.

If $|\{j_{1},\cdots,j_{d}\}|=1$, we have $|G_{T_{g^{s}}^{j_{1},\cdots,j_{d}}}|=o\left(g\right)=k$
and the cardinality of the orbit of $T_{g^{s}}^{j_{1},\cdots,j_{d}}$
is
\[
\left|\mathcal{O}_{T_{g^{s}}^{j_{1},\cdots,j_{d}}}\right|=\left[G:G_{T_{g^{s}}^{j_{1},\cdots,j_{d}}}\right]=1.
\]
In this case we have $\Lambda_{G_{T_{g^{s}}^{j_{1},\cdots,j_{d}}}}=\left\{ \lambda_{n}\mid n=0,\cdots,k-1\right\} $
where $\lambda_{n}\left(g\right)=e^{-\frac{2\pi in}{k}}$.

Let $I$ be a subset of $\left[0,p\right]^{k}\backslash\left\{ \left(i,\cdots,i\right)\mid i\in\left[0,p\right]\right\} $
consisting of the orbit representatives under the action of $G=\left\langle g\right\rangle .$
Then it is easy to see the following proposition.

\begin{proposition} The irreducible $\left(V^{\otimes k}\right)^{\left\langle g\right\rangle }$-modules
consist of
\begin{align*}
 & \left\{ W^{i_{1},\cdots,i_{k}}\mid\left(i_{1},\cdots,i_{k}\right)\in I\right\} ,\\
 & \left\{ \left(W^{j,\cdots,j}\right)^{n}\mid j\in\left[0,p\right],0\le n<k\right\} ,\\
 & \left\{ \left(T_{g^{s}}^{j_{1},\cdots,j_{d}}\right)^{n}\mid1\le s<k,j_{1},\cdots,j_{d}\in\left[0,p\right],|\{j_{1},\cdots,j_{d}\}|=1,0\le n<k,d=\text{gcd}\left(s,k\right),k=dl\right\} ,\\
 & \left\{ \left(T_{g^{s}}^{j_{1},\cdots,j_{d}}\right)^{t}\mid1\le s<k,j_{1},\cdots,j_{d}\in\left[0,p\right],j_{1}\le\cdots\le j_{d},|\{j_{1},\cdots,j_{d}\}|>1,0\le t<l,d=\text{gcd}\left(s,k\right),k=dl\right\} .
\end{align*}

\end{proposition}

\subsection{$S$-matrix of $\left(V^{\otimes k}\right)^{\left\langle g\right\rangle }$}

Let $J$ be as before and $i,j\in J$. Let $M^{i}$ be a $g_{i}$-twisted
$V^{\otimes k}$-module and $M^{j}$ be a $g_{j}$-twisted $V^{\otimes k}$-module,
where $g_{i},g_{j}\in\left\langle g\right\rangle $. Let $C_{i,j}$
be the least subset of $\left\langle g\right\rangle $ such that
\[
\left\{ M^{j}\circ k\mid k\in C_{i,j}\right\} =\mathcal{O}_{j}\cap\left(\cup_{h\in G_{M^{i}}}\mathcal{\mathfrak{U}}\left(h,g_{i}^{-1}\right)\right).
\]
Note that here $G=\left\langle g\right\rangle $ is an abelian group.
We will use the following result which is from Corollary 5.4 of \cite{DRX21}:

\begin{lemma} \label{orbifold S-matrix formula} Let $i,j$ be as
before, $\lambda\in\Lambda_{G_{M^{i}}}$ and $\mu\in\Lambda_{G_{M^{j}}}$.
Then
\begin{equation}
S_{M_{\lambda}^{i},M_{\mu}^{j}}=\frac{1}{\left|G_{M^{i}}\right|}\sum_{k\in C_{i,j}}S_{M^{i},M^{j}\circ k}\overline{\lambda\left(\overline{g_{j}}\right)}\mu\left(\overline{g_{i}^{-1}}\right)\label{Restricted S-matrix formula}
\end{equation}
if $C_{i,j}$ is not empty, and $S_{M_{\lambda}^{i},M_{\mu}^{j}}=0$
otherwise, where $\overline{g_{j}}$ is the element in twisted group
algebra $\mathbb{C}^{\alpha_{\mathcal{M}}}\left[G_{M}\right]$ corresponding
to $g_{j}\in G_{M}$ and $\bar{x}$ is the complex conjugate of the
complex number $x$.

\end{lemma}

Now we give explicit expressions of the entries of the $S$-matrix
of $\left(V^{\otimes k}\right)^{\left\langle g\right\rangle }$, where
we use the fact that $S$-matrix is symmetric \cite{H}. The entries
of the $S$-matrix of $\left(V^{\otimes k}\right)^{\left\langle g\right\rangle }$
that only involve untwisted $V^{\otimes k}$-modules (i.e., (3) and
(4) in the following theorem) are given in Theorem 6.8 of \cite{DRX21}.  For completeness we list them here.

\begin{theorem} Let $1\le r,s<k$, $d=\text{gcd}\left(s,k\right),$
$d_{1}=\text{gcd}\left(r,k\right)$ and $f=\text{gcd}$$\left(d,r\right)$.
Set $l=\frac{k}{d},$ $l_{1}=\frac{k}{d_{1}}$, $b=\frac{d}{f}$,
$a=\frac{d_{1}}{f}.$ Suppose $i_{1},\cdots,i_{f},j_{1},\cdots,j_{f}\in\left[0,p\right],$
$j_{1}\le\cdots\le j_{f}$.


(1) If $|\{i_{1},\cdots,i_{f}\}|>1$ and $i_{1}\le\cdots\le i_{f}$,
then for $0\le t<l_{1}$, $0\le t_{1}<l$, $0\le n<k$,
\[
S_{\left(T_{g^{r}}^{i_{1},\cdots,i_{f};a}\right)^{t},N}=\frac{1}{l_{1}}\begin{cases}
e^{\frac{2\pi i\left(st+rt_{1}\right)}{k}}\sum_{i=0}^{d-1}S_{T_{g^{r}}^{i_{1},\cdots,i_{f};a},T_{g^{s}}^{j_{1},\cdots,j_{f};b}\circ g^{i}}, & \text{if}\ N=\left(T_{g^{s}}^{j_{1},\cdots,j_{f};b}\right)^{t_{1}},\left|\left\{ j_{1},\cdots,j_{f}\right\} \right|>1\\
 & \ \ \ \ \text{and\ \ }d_{1}m\equiv s\  (\text{mod}\  k)\ \text{for\ some}\ m;\\
e^{\frac{2\pi i\left(st+rn\right)}{k}}S_{T_{g^{r}}^{i_{1},\cdots,i_{f};a},T_{g^{s}}^{j_{1},\cdots,j_{f};b},} & \text{if}\ N=\left(T_{g^{s}}^{j_{1},\cdots,j_{f};b}\right)^{n},\left|\left\{ j_{1},\cdots,j_{f}\right\} \right|=1\\
 & \ \ \ \ \text{and\ }d_{1}m\equiv s\  (\text{mod}\ k)\ \text{for\ some}\ m;\\
\sum_{i=0}^{r-1}S_{T_{g^{r}}^{i_{1},\cdots,i_{f};a},\left(\left(W^{j_{1}}\right)^{\otimes l_{1}}\otimes\cdots\otimes\left(W^{j_{d_{1}}}\right)^{\otimes l_{1}}\right)\circ g^{i}}, & \text{if}\ N=\left(W^{j_{1}}\right)^{\otimes l_{1}}\otimes\cdots\otimes\left(W^{j_{d_{1}}}\right)^{\otimes l_{1}}\\
 & \ \ \ \ \text{and\ \ }d_{1}>1;\\
e^{\frac{2\pi irn}{k}}S_{T_{g^{r}}^{W^{i_{1}}},W^{j,\cdots,j}}, & \text{if}\ N=\left(W^{j,\cdots,j}\right)^{n}\ \text{and}\ d_{1}=f=1;\\
0 & \text{otherwise},
\end{cases}
\]
where $S_{T_{g^{r}}^{i_{1},\cdots,i_{f};a},T_{g^{s}}^{j_{1},\cdots,j_{f};b}}$,
$S_{T_{g^{r}}^{i_{1},\cdots,i_{f};a},\left(W^{j_{1}}\right)^{\otimes l_{1}}\otimes\cdots\otimes\left(W^{j_{d_{1}}}\right)^{\otimes l_{1}}}$
and $S_{T_{g^{r}}^{W^{i_{1}}},W^{j,\cdots,j}}$ are given in Lemmas
\ref{S-matrix: twisted and twisted} and \ref{S-matrix, nontwisted and twisted}.

(2) If $|\{i_{1},\cdots,i_{f}\}|=1$, then for $0\le t< k,$ $0\le t_{1}<l$,  $0\le n<k$,
\[
S_{\left(T_{g^{r}}^{i_{1},\cdots,i_{f};a}\right)^{t},N}=\frac{1}{k}\begin{cases}
e^{\frac{2\pi i\left(st+rt_{1}\right)}{k}}\sum_{i=0}^{d-1}S_{T_{g^{r}}^{i_{1},\cdots,i_{f};a},T_{g^{s}}^{j_{1},\cdots,j_{f};b}\circ g^{i}}, & \text{if}\ N=\left(T_{g^{s}}^{j_{1},\cdots,j_{f};b}\right)^{t_{1}},\left|\left\{ j_{1},\cdots,j_{f}\right\} \right|>1;\\
e^{\frac{2\pi i\left(st+rn\right)}{k}}S_{T_{g^{r}}^{i_{1},\cdots,i_{f};a},T_{g^{s}}^{j_{1},\cdots,j_{f};b},} & \text{if}\ N=\left(T_{g^{s}}^{j_{1},\cdots,j_{f};b}\right)^{n},\left|\left\{ j_{1},\cdots,j_{f}\right\} \right|=1;\\
\sum_{i=0}^{r-1}S_{T_{g^{r}}^{i_{1},\cdots,i_{f};a},\left(\left(W^{j_{1}}\right)^{\otimes l_{1}}\otimes\cdots\otimes\left(W^{j_{d_{1}}}\right)^{\otimes l_{1}}\right)\circ g^{i}}, & \text{if}\ N=\left(W^{j_{1}}\right)^{\otimes l_{1}}\otimes\cdots\otimes\left(W^{j_{d_{1}}}\right)^{\otimes l_{1}}\\
 & \ \ \ \ \text{and\ \ }d_{1}>1;\\
e^{\frac{2\pi irn}{k}}S_{T_{g^{r}}^{W^{i_{1}}},W^{j,\cdots,j}}, & \text{if}\ N=\left(W^{j,\cdots,j}\right)^{n}\ \text{and}\ d_{1}=f=1;\\
0 & \text{otherwise},
\end{cases}
\]
where $S_{T_{g^{r}}^{i_{1},\cdots,i_{f};a},T_{g^{s}}^{j_{1},\cdots,j_{f};b}}$,
$S_{T_{g^{r}}^{i_{1},\cdots,i_{f};a},\left(W^{j_{1}}\right)^{\otimes l_{1}}\otimes\cdots\otimes\left(W^{j_{d_{1}}}\right)^{\otimes l_{1}}}$
and $S_{T_{g^{r}}^{W^{i_{1}}},W^{j,\cdots,j}}$ are given in Lemmas
\ref{S-matrix: twisted and twisted} and \ref{S-matrix, nontwisted and twisted}.

(3) Let $\left(i_{1},\cdots,i_{k}\right),$$\left(j_{1},\cdots,j_{k}\right)\in I$,
$j\in\left[0,p\right]$ and $0\le n<k$. Then
\[
S_{W^{i_{1},\cdots,i_{k}},N}=\begin{cases}
\sum_{n=0}^{k-1}\prod_{t=1}^{k}S_{W^{i_{t}}, W^{j_{t}+n}}, & \ \text{if}\ N=W^{j_{1},\cdots,j_{k}};\\
\prod_{t=1}^{k}S_{W^{i_{t}},W^{j}}, & \text{if}\ N=\left(W^{j,\cdots,j}\right)^{n}.
\end{cases}
\]

(4) Let $i,j\in\left[0,p\right]$ and $0\le m,n<k$. We have
\[
S_{\left(W^{i,\cdots,i}\right)^{m},\left(W^{j,\cdots,j}\right)^{n}}=\frac{1}{k}S_{W^{i}, W^{j}}^{k}.
\]
\end{theorem} \begin{proof} (1)  We will use notations from (\ref{Restricted S-matrix formula}). Take     $M^{i}=T_{g^{r}}^{i_{1},\cdots,i_{f};a}$ and
$g_{i}=g^{r}$. Then $G_{M^{i}}=\left\langle g^{d_{1}}\right\rangle $
where $d_{1}=\text{gcd}\left(r,k\right)$ and $\left|G_{M^{i}}\right|=l_{1}$.
Now we have $\Lambda_{G_{M^{i}}}=\left\{ \mu_{t}\mid0\le t<l_{1}\right\} $
where $\mu_{t}\left(g^{d_{1}}\right)=e^{-\frac{2\pi it}{l_{1}}}.$
By Lemma \ref{orbifold S-matrix formula}, it suffices to find $\overline{\lambda\left(\overline{g_{j}}\right)}\mu\left(\overline{g_{i}^{-1}}\right)$
and $C_{ij}$.
Now we have
\begin{equation}
\cup_{h\in G_{M^{i}}}\mathcal{\mathfrak{U}}\left(h,g_{i}^{-1}\right)=\cup_{h\in\left\langle g^{d_{1}}\right\rangle }\mathcal{\mathfrak{U}}\left(h,g^{-r}\right)=\cup_{m\in\mathbb{N}}\mathcal{\mathfrak{U}}\left(g^{d_{1}m},g^{-r}\right).\label{eq:U}
\end{equation}

{\it Case 1.} If $N=\left(T_{g^{s}}^{j_{1},\cdots,j_{f};b}\right)^{t_{1}}$ and
$\left|\left\{ j_{1},\cdots,j_{f}\right\} \right|>1,$ by taking $M^{j}=T_{g^{s}}^{j_{1},\cdots,j_{f};b},$
$g_{j}=g^{s}$, $\lambda=\lambda_{t}$ and $\mu=\mu_{t_{1}}$, we
obtain $\overline{\lambda_{t}\left(g^{s}\right)}\mu_{t_{1}}\left(\overline{g^{-r}}\right)=e^{\frac{2\pi i\left(st+rt_{1}\right)}{k}}$
and
\[
\mathcal{O}_{j}\cap\left(\cup_{h\in G_{M^{i}}}\mathcal{\mathfrak{U}}\left(h,g_{i}^{-1}\right)\right)=\left\{ T_{g^{s}}^{j_{1},\cdots,j_{f};b}\circ g^{t}\mid0\le t<d\right\} \cap\left(\cup_{n\in\mathbb{N}}\mathcal{\mathfrak{U}}\left(g^{d_{1}n},g^{-r}\right)\right).
\]
By Lemma \ref{g^s, g^r modules}, $T_{g^{s}}^{j_{1},\cdots,j_{f};b}\circ g^{t}\in\mathcal{\mathfrak{U}}\left(g^{s},g^{-r}\right)$
for all $0\le t<d$. If there exists $m\in\mathbb{N}$ such that $d_{1}m\equiv s$
(mod $k$), then $\left|C_{ij}\right|=\left|\mathcal{O}_{j}\right|=\left|M^{j}\circ G\right|=d$.
Using (\ref{Restricted S-matrix formula}), we obtain
\[
S_{\left(T_{g^{r}}^{i_{1},\cdots,i_{f};a}\right)^{t},\left(T_{g^{s}}^{j_{1},\cdots,j_{f};b}\right)^{t_{1}}}=\frac{1}{l_{1}}e^{\frac{2\pi i\left(st+rt_{1}\right)}{k}}\sum_{i=0}^{d-1}S_{T_{g^{r}}^{i_{1},\cdots,i_{f};a},T_{g^{s}}^{j_{1},\cdots,j_{f};b}\circ g^{i}}.
\]
If $d_{1}m\not\equiv s$ (mod $k$) for any $m,$ then $C_{ij}$
is an empty set and the corresponding entry in the $S$-matrix is
$0$.

{\it Case 2.} If $N=\left(T_{g^{s}}^{j_{1},\cdots,j_{f};b}\right)^{n}$ and $\left|\left\{ j_{1},\cdots,j_{f}\right\} \right|=1,$
we take $M^{j}=T_{g^{s}}^{j_{1},\cdots,j_{f};b}$ and $g_{j}=g^{s}$.
Then it is clear that $T_{g^{s}}^{j_{1},\cdots,j_{f};b}\circ g^{i}\cong T_{g^{s}}^{j_{1},\cdots,j_{f};b}$
for any $0\le i<k$. Thus $\left|\mathcal{O}_{j}\right|=1$. By taking
$\lambda=\lambda_{t}$ and $\mu=\mu_{n}$, we obtain $\overline{\lambda_{t}\left(g^{s}\right)}\mu_{n}\left(\overline{g^{-r}}\right)=e^{\frac{2\pi i\left(st+rn\right)}{k}}$.
By Lemma \ref{g^s, g^r modules}, $T_{g^{s}}^{j_{1},\cdots,j_{f};b}\circ g^{i}\in\mathcal{\mathfrak{U}}\left(g^{s},g^{-r}\right)$
for all $0\le i<d$. If there exists $m\in\mathbb{N}$ such that $d_{1}m\equiv s$
(mod $k$). Then $\left|C_{ij}\right|=\left|\mathcal{O}_{j}\right|=1$. By (\ref{Restricted S-matrix formula}), we get
\[
S_{\left(T_{g^{r}}^{i_{1},\cdots,i_{f};a}\right)^{t},\left(T_{g^{s}}^{j_{1},\cdots,j_{f};b}\right)^{n}}=\frac{1}{l_{1}}e^{\frac{2\pi i\left(st+rn\right)}{k}}S_{T_{g^{r}}^{i_{1},\cdots,i_{f};a},T_{g^{s}}^{j_{1},\cdots,j_{f};b}.}
\]
If $d_{1}m\not\equiv s$ (mod $k$) for any $m,$ then $C_{ij}$
is an empty set and the corresponding entry in the $S$-matrix is
$0$.

{\it Case 3.} If $N=\left(W^{j_{1}}\right)^{\otimes l_{1}}\otimes\cdots\otimes\left(W^{j_{d_{1}}}\right)^{\otimes l_{1}}$
with $d_{1}>1$, we take $M^{j}=\left(W^{j_{1}}\right)^{\otimes l_{1}}\otimes\cdots\otimes\left(W^{j_{d_{1}}}\right)^{\otimes l_{1}}$
and $g_{j}=1$. Recall that $M^{j}\in\mathfrak{U}\left(1,g^{r}\right).$
Now we have
\[
\mathcal{O}_{j}\cap\left(\cup_{h\in G_{M^{i}}}\mathcal{\mathfrak{U}}\left(h,g_{i}^{-1}\right)\right)=\left\{ \left(\left(W^{j_{1}}\right)^{\otimes l_{1}}\otimes\cdots\otimes\left(W^{j_{d_{1}}}\right)^{\otimes l_{1}}\right)\circ g^{t}\mid0\le t<r\right\} \cap\left(\cup_{n\in\mathbb{N}}\mathcal{\mathfrak{U}}\left(g^{d_{1}n},g^{-r}\right)\right).
\]
By definition of $d_{1}$, it is clear that there exists $m\in\mathbb{N}$
such that $d_{1}m\equiv 0$ (mod $k$). By (\ref{1,g^s}), for any $0\le i<r$,   $\left(\left(W^{j_{1}}\right)^{\otimes l_{1}}\otimes\cdots\otimes\left(W^{j_{d_{1}}}\right)^{\otimes l_{1}}\right)\circ g^{i}\in\mathfrak{U}\left(1,g^{-r}\right)$.  Therefore, $\left|C_{ij}\right|=\left|\mathcal{O}_{j}\right|=r$.
By taking $\lambda=\lambda_{t}$ and $\mu=1$, we get $\overline{\lambda_{t}\left(1\right)}\mu\left(\overline{g^{-r}}\right)=1$.
Using (\ref{Restricted S-matrix formula}), we obtain
\[
S_{\left(T_{g^{r}}^{i_{1},\cdots,i_{f};a}\right)^{t},\left(W^{j_{1}}\right)^{\otimes l_{1}}\otimes\cdots\otimes\left(W^{j_{d_{1}}}\right)^{\otimes l_{1}}}=\frac{1}{l_{1}}\sum_{i=0}^{r-1}S_{T_{g^{r}}^{i_{1},\cdots,i_{f};a}, \left(\left(W^{j_{1}}\right)^{\otimes l_{1}}\otimes\cdots\otimes\left(W^{j_{d_{1}}}\right)^{\otimes l_{1}}\right)\circ g^{i}}.
\]

{\it Case 4.} If $N=\left(W^{j,\cdots,j}\right)^{n},$ we take $M^{j}=W^{j,\cdots,j},$
$g_{j}=1.$ Note that $\left|\mathcal{O}_{j}\right|=1$ since $M^{j}\circ g^{i}\cong M^{j}$
for any $0\le i<k$. Similar to the previous case, there exists $m\in\mathbb{N}$
such that $d_{1}m\equiv 0$ (mod $k$) and it is obvious that each $W^{j,\cdots,j}\in\mathcal{\mathfrak{U}}\left(1,g^{-r}\right)$.
So now we have $\left|C_{ij}\right|=\left|\mathcal{O}_{j}\right|=1$.
By taking $\lambda=\lambda_{t}$ and $\mu=\mu_{n}$, we get $\overline{\lambda_{t}\left(1\right)}\mu_{n}\left(\overline{g^{-r}}\right)=e^{\frac{2\pi irn}{k}}$.
Notice that if $d_{1}=f=1$, $T_{g^{r}}^{i_{1},\cdots,i_{f};a}$
can be written as $T_{g^{r}}^{W^{i_{1}}}.$ Thus we obtain
\[
S_{\left(T_{g^{r}}^{i_{1},\cdots,i_{f};a}\right)^{t},\left(W^{j,\cdots,j}\right)^{n}}=\frac{1}{l_{1}}e^{\frac{2\pi irn}{k}}S_{T_{g^{r}}^{W^{i_{1}}},W^{j,\cdots,j}}.
\]

{\it Case 5.} If the $g^{s}$-twisted $V^{\otimes k}$-module $M^{j}$ given in (\ref{g^s-module general form})  is not of the form $T_{g^{s}}^{j_{1},\cdots,j_{f};b}$, that is,
if the $M^{j}$ is not $g^{r}$-stable, then
\[
\mathcal{O}_{j}\cap\left(\cup_{h\in G_{M^{i}}}\mathcal{\mathfrak{U}}\left(h,g_{i}^{-1}\right)\right)=\left\{ M^{j}\circ g^{t}\mid0\le t<k\right\} \cap\left(\cup_{n\in\mathbb{N}}\mathcal{\mathfrak{U}}\left(g^{d_{1}n},g^{-r}\right)\right)=\emptyset
\]
and hence $C_{ij}=\emptyset$. Thus we have $S_{\left(T_{g^{r}}^{i_{1},\cdots,i_{f};a}\right)^{t},\left(M^{j}\right)^{t_{1}}}=0$
 for any $0\le t_{1}<l$, where $(M^j)^{t_1}$ is given in (\ref{g^s-decomposition}). Similarly, if $M^{j}=W^{j_{1}}\otimes\cdots\otimes W^{j_{k}}$
cannot be written in the form $\left(W^{j_{1}}\right)^{\otimes l_{1}}\otimes\cdots\otimes\left(W^{j_{d_{1}}}\right)^{\otimes l_{1}}$
or $M^{i,\cdots,i}$ for $\left\{ j_{1},\cdots,j_{d_{1}}\right\} \subset\left\{ 0,\cdots,p\right\} $,
$i\in\left\{ 0,\cdots,p\right\} $ , then $M^{j}$ is not $g^{r}$
-stable and hence we also have $C_{ij}=0$ for this case. Therefore
the corresponding entries in the $S$-matrix are zeros.

(2) We will use  Lemma \ref{orbifold S-matrix formula}. Take $M^{i}=T_{g^{r}}^{i_{1},\cdots,i_{f};a}$ with $|\{i_{1},\cdots,i_{f}\}|=1$
and $g_{i}=g^{r}$ in (\ref{Restricted S-matrix formula}). Then $G_{M^{i}}=\left\langle g\right\rangle $
and hence $\left|G_{M^{i}}\right|=k$. It suffices to find $\overline{\lambda\left(\overline{g_{j}}\right)}\mu\left(\overline{g_{i}^{-1}}\right)$
and $C_{ij}$. Note that $\Lambda_{G_{T_{g^{r}}^{i_{1},\cdots,i_{f};a}}}=\left\{ \mu_{n}\mid n=0,\cdots,k-1\right\} $
where $\mu_{n}\left(g\right)=e^{-\frac{2\pi in}{k}}$. In  this
case we have
\[
\cup_{h\in G_{M^{i}}}\mathcal{\mathfrak{U}}\left(h,g_{i}^{-1}\right)=\cup_{h\in\left\langle g\right\rangle }\mathcal{\mathfrak{U}}\left(h,g^{-r}\right)=\cup_{n\in\mathbb{N}}\mathcal{\mathfrak{U}}\left(g^{n},g^{-r}\right).
\]

{\it Case 1. } If $N=\left(T_{g^{s}}^{j_{1},\cdots,j_{f};b}\right)^{t_{1}}$ and
$\left|\left\{ j_{1},\cdots,j_{f}\right\} \right|>1,$ by taking $M^{j}=T_{g^{s}}^{j_{1},\cdots,j_{f};b}$,
$g_{j}=g^{s}$, $\lambda=\lambda_{t}$ and $\mu=\mu_{t_{1}}$, we
obtain $\overline{\lambda_{t}\left(g^{s}\right)}\mu_{t_{1}}\left(\overline{g^{-r}}\right)=e^{\frac{2\pi i\left(st+rt_{1}\right)}{k}}$
and
\begin{align*}
 & \mathcal{O}_{j}\cap\left(\cup_{h\in G_{M^{i}}}\mathcal{\mathfrak{U}}\left(h,g_{i}^{-1}\right)\right)\\
= & \left\{ T_{g^{s}}^{j_{1},\cdots,j_{f};b}\circ g^{t}\mid0\le t<d\right\} \cap\left(\cup_{n\in\mathbb{N}}\mathcal{\mathfrak{U}}\left(g^{n},g^{-r}\right)\right)=\left\{ T_{g^{s}}^{j_{1},\cdots,j_{f};b}\circ g^{t}\mid0\le t<d\right\} =\mathcal{O}_{j}.
\end{align*}
Hence $\left|C_{ij}\right|=\left|\mathcal{O}_{j}\right|=\left|M^{j}\circ G\right|=d$.
Using (\ref{Restricted S-matrix formula}), we obtain
\[
S_{\left(T_{g^{r}}^{i_{1},\cdots,i_{f};a}\right)^{t},\left(T_{g^{s}}^{j_{1},\cdots,j_{f};b}\right)^{t_{1}}}=\frac{1}{k}e^{\frac{2\pi i\left(st+rt_{1}\right)}{k}}\sum_{i=0}^{d-1}S_{T_{g^{r}}^{i_{1},\cdots,i_{f};a},T_{g^{s}}^{j_{1},\cdots,j_{f};b}\circ g^{i}}.
\]

{\it Case 2.} If $N=\left(T_{g^{s}}^{j_{1},\cdots,j_{f};b}\right)^{n}$ and $\left|\left\{ j_{1},\cdots,j_{f}\right\} \right|=1,$
by taking $M^{j}=T_{g^{s}}^{j_{1},\cdots,j_{f};b},$ $g_{j}=g^{s}$,
$\lambda=\lambda_{t}$ and $\mu=\mu_{n}$, we obtain $\overline{\lambda_{t}\left(g^{s}\right)}\mu_{n}\left(\overline{g^{-r}}\right)=e^{\frac{2\pi i\left(st+rn\right)}{k}}$.
Since in this case $\left|C_{ij}\right|=\left|\mathcal{O}_{j}\right|=1$,
we get
\[
S_{\left(T_{g^{r}}^{i_{1},\cdots,i_{f};a}\right)^{t},\left(T_{g^{s}}^{j_{1},\cdots,j_{f};b}\right)^{n}}=\frac{1}{k}e^{\frac{2\pi i\left(st+rn\right)}{k}}S_{T_{g^{r}}^{i_{1},\cdots,i_{f};a},T_{g^{s}}^{j_{1},\cdots,j_{f};b}.}
\]

{\it Case 3.} If $N=\left(W^{j_{1}}\right)^{\otimes l_{1}}\otimes\cdots\otimes\left(W^{j_{d_{1}}}\right)^{\otimes l_{1}}$
with $d_{1}>1$, we take $M^{j}=\left(W^{j_{1}}\right)^{\otimes l}\otimes\cdots\otimes\left(W^{j_{d_{1}}}\right)^{\otimes l}$
and $g_{j}=1$. By (\ref{1,g^s}),  $\left(\left(W^{j_{1}}\right)^{\otimes l_{1}}\otimes\cdots\otimes\left(W^{j_{d_{1}}}\right)^{\otimes l_{1}}\right)\circ g^{i}\in\mathfrak{U}\left(1,g^{-r}\right)$
for any $0\le i<r$. Therefore,
\[
\mathcal{O}_{j}\cap\left(\cup_{h\in G_{M^{i}}}\mathcal{\mathfrak{U}}\left(h,g_{i}^{-1}\right)\right)=\left\{ \left(\left(W^{j_{1}}\right)^{\otimes l_{1}}\otimes\cdots\otimes\left(W^{j_{d_{1}}}\right)^{\otimes l_{1}}\right)\circ g^{i}\mid0\le i<r\right\} \cap\left(\cup_{n\in\mathbb{N}}\mathcal{\mathfrak{U}}\left(g^{n},g^{-r}\right)\right)=\mathcal{O}_{j}
\]
and $\left|C_{ij}\right|=r$. By taking $\lambda=\lambda_{t}$ and
$\mu=1$, we get $\overline{\lambda\left(\overline{g_{j}}\right)}\mu\left(\overline{g_{i}^{-1}}\right)=1$
and hence
\[
S_{\left(T_{g^{r}}^{i_{1},\cdots,i_{f};a}\right)^{t},\left(W^{j_{1}}\right)^{\otimes l_{1}}\otimes\cdots\otimes\left(W^{j_{d_{1}}}\right)^{\otimes l_{1}}}=\frac{1}{k}\sum_{i=0}^{r-1}S_{T_{g^{r}}^{i_{1},\cdots,i_{f};a},\left(\left(W^{j_{1}}\right)^{\otimes l_{1}}\otimes\cdots\otimes\left(W^{j_{d_{1}}}\right)^{\otimes l_{1}}\right)\circ g^{i}}.
\]

{\it Case 4.} If $N=\left(W^{j,\cdots,j}\right)^{n},$ we take $M^{j}=W^{j,\cdots,j}$
and $g_{j}=1.$ Note that $\left|\mathcal{O}_{j}\right|=1$ since
$W^{j,\cdots,j}$ is $g^{i}$-stable for any $0\le i<k$. It is clear
that $W^{j,\cdots,j}\circ g^{i}\in\mathfrak{U}\left(1,g^{-r}\right)$
for any $0\le i<k$. Thus
\[
\mathcal{O}_{j}\cap\left(\cup_{h\in G_{M^{i}}}\mathcal{\mathfrak{U}}\left(h,g_{i}^{-1}\right)\right)=\left\{ W^{j,\cdots,j}\circ g^{i}\mid0\le i<k\right\} \cap\left(\cup_{n\in\mathbb{N}}\mathcal{\mathfrak{U}}\left(g^{n},g^{-r}\right)\right)=\mathcal{O}_{j}
\]
and $\left|C_{ij}\right|=\left|\mathcal{O}_{j}\right|=1$. By taking
$\lambda=\lambda_{t}$ and $\mu=\mu_{n}$, we obtain $\overline{\lambda_{t}\left(1\right)}\mu_{n}\left(\overline{g^{-r}}\right)=e^{\frac{2\pi irn}{k}}$.
Also note that when $d_{1}=f=1$, $T_{g^{r}}^{i_{1},\cdots,i_{f};a}$
can be written in  the form $T_{g^{r}}^{W^{i_{1}}}.$ Now we get

\[
S_{\left(T_{g^{r}}^{i_{1},\cdots,i_{f};a}\right)^{t},\left(W^{j,\cdots,j}\right)^{n}}=\frac{1}{k}e^{\frac{2\pi irn}{k}}S_{T_{g^{r}}^{W^{i_{1}}},W^{j,\cdots,j}}.
\]

{\it Case 5.} The proof for other cases are similar to the proof for {\it Case 5} in (1). \end{proof}

\vskip10pt {\footnotesize{}{ }\textbf{\footnotesize{}C. Dong}{\footnotesize{}:
Department of Mathematics, University of California Santa Cruz, CA 95064 USA; }\texttt{dong@ucsc.edu}{\footnotesize \par}

\textbf{\footnotesize{}F. Xu}{\footnotesize{}:  Department of Mathematics, University of California, Riverside, CA 92521 USA; }\texttt{xufeng@math.ucr.edu}{\footnotesize \par}

\textbf{\footnotesize{}N. Yu}{\footnotesize{}: School of Mathematical
Sciences, Xiamen University, Fujian, 361005, CHINA;} \texttt{
ninayu@xmu.edu.cn}{\footnotesize \par}
\end{document}